\documentclass[twoside,a4paper,reqno,11pt]{amsart}
\usepackage{amsfonts, amsbsy, amsmath, amsthm, amssymb, latexsym, verbatim, enumerate}
\usepackage{mathrsfs}
\usepackage[top=30mm,right=30mm,bottom=30mm,left=30mm]{geometry}

\usepackage{bm}
\usepackage[pdftex]{color,graphicx}

\makeatletter
\newcommand{\imod}[1]{\allowbreak\mkern4mu({\operator@font mod}\,\,#1)}
\makeatother

\def\hal{\unskip\nobreak\hfill\penalty50\hskip10pt\hbox{}\nobreak

\hfill\vrule height 5pt width 6pt depth 1pt\par\vskip 2mm}

\usepackage{bbm}
\usepackage[all,cmtip]{xy}

\usepackage{enumerate}

 \def\Z{\mathbb Z}
  \def\N{\mathbb N}
  \def\F{\mathbb F}
  \def\R{\mathbb R}

  \def\C{\mathbb C}

\def\c{\chi}
  \def\a{\alpha}
  \def\b{\beta}
  \def\g{\gamma}
  \def\G{\Gamma}

  \def\e{\epsilon}
  \def\la{\lambda}

  \def\no{\noindent}
  
  \def\pf{\noindent {\bf Proof $\;$ }}

\newcommand{\IBR}{{\mathrm {IBr}}}

\newcommand{\diag}{{\mathrm {diag}}}

\newcommand{\Irr}{{\mathrm {Irr}}}

\newcommand{\Deg}{{\mathrm {Deg}}}

\newcommand{\tw}[1]{{}^#1\!}
\newcommand{\SR}{\tw* R}
\newcommand{\SSS}{{\sf S}}

\newcommand{\St}{{\sf {St}}}

\newcommand{\CB}{{\bf C}}
\newcommand{\ZB}{{\bf Z}}
\newcommand{\NB}{{\bf N}}

\newcommand{\KK}{{\mathbb K}}

\newcommand{\GC}{{\mathcal G}}
\newcommand{\PC}{{\mathcal P}}
\newcommand{\TC}{{\mathcal T}}

\newcommand{\ZC}{{\mathcal Z}}

\newcommand{\TG}{\tilde\Gamma}

\newcommand{\HC}{{\mathcal H}}
\newcommand{\UC}{{\mathcal U}}

\newcommand{\LC}{{\mathcal L}}

\newcommand{\OC}{{\mathcal O}}
\newcommand{\OCS}{{\mathcal O}^*}

\newcommand{\GCF}{\GC^{F}}

\newcommand{\al}{\alpha}

\newcommand{\gam}{\gamma}
\newcommand{\lam}{\lambda}

\newcommand{\tchi}{\tilde\chi}
\newcommand{\supp}{{\mathrm {supp}}}

\numberwithin{equation}{section}

  \def\hal{\unskip\nobreak\hfil\penalty50\hskip10pt\hbox{}\nobreak
  \hfill\vrule height 5pt width 6pt depth 1pt\par\vskip 2mm}

\begin{document}

\author[R. Bezrukavnikov]{Roman Bezrukavnikov}
\address{R. Bezrukavnikov,   Department of Mathematics,     MIT,     Cambridge, MA 02139, USA}
\email{bezrukav@math.mit.edu}

\author[M.W. Liebeck]{Martin W. Liebeck}
\address{M.W. Liebeck, Department of Mathematics,
    Imperial College, London SW7 2BZ, UK}
\email{m.liebeck@imperial.ac.uk}

\author[A. Shalev]{Aner Shalev}
\address{A. Shalev, Institute of Mathematics, Hebrew University, Jerusalem 91904, Israel}
\email{shalev@math.huji.ac.il}

\author[P. H. Tiep]{Pham Huu Tiep}
\address{P.H. Tiep, Department of Mathematics,    University of Arizona,    Tucson, AZ 85721, USA}
\email{tiep@math.arizona.edu}

\title[Character bounds for finite groups of Lie type]{Character bounds for finite groups of Lie type}



  \maketitle

\begin{abstract}
We establish new bounds on character values and character ratios for finite groups $G$ of Lie type, which are considerably stronger than previously known bounds, and which are best possible
in many cases. These bounds have the form $|\chi(g)| \le c\chi(1)^{\a_g}$, and give rise to a variety
of applications, for example to covering numbers and mixing times of random walks on such groups.
In particular we deduce that, if $G$ is a classical group in dimension $n$, then, under some conditions
on $G$ and $g \in G$, the mixing time of the random walk on $G$ with the conjugacy class of $g$ as a
generating set is (up to a small multiplicative constant) $n/s$, where $s$ is the support of $g$.
\end{abstract}

\footnotetext{The first author was partially supported by the NSF grants DMS-1102434 and DMS-1601953.}
\footnotetext{The second and third authors acknowledge the support of EPSRC grant EP/H018891/1.}
\footnotetext{The third author acknowledges the support of ERC advanced grant
247034, ISF grant 1117/13 and the Vinik chair of mathematics which he holds.}
\
\footnotetext{The fourth author was partially supported by the NSF grants DMS-1201374 and DMS-1665014, the Simons Foundation
Fellowship 305247, the EPSRC, and the Mathematisches Forschungsinstitut Oberwolfach.
Parts of the paper were written while the fourth author visited the Departments of Mathematics of Imperial College, London, and
University of Chicago. It is a pleasure to thank Imperial College, Prof. Ngo Bao Chau, and the University of Chicago
for generous hospitality and stimulating environment.}
\footnotetext{The authors thank M. Geck, G. Malle, J. Michel, and especially J. Taylor for many helpful discussions.}

  \newtheorem{thm}{Theorem}[section]
  \newtheorem{prop}[thm]{Proposition}
  \newtheorem{lem}[thm]{Lemma}
  \newtheorem{cor}[thm]{Corollary}
   \newtheorem{rem}[thm]{Remark}
   \newtheorem{exa}[thm]{Example}

\tableofcontents


\section{Introduction}

For a finite group $G$, a {\it character ratio} is a complex number of the form $\frac{\c(g)}{\c(1)}$, where $g\in G$ and $\c$ is an irreducible character of $G$. Upper bounds for absolute values of character values and character ratios have long been of interest, for various reasons; these include applications to random generation, covering numbers,
mixing times of random walks, the study of word maps, representation varieties and other areas.
For example, character ratios are connected with the well-known formula
\[
\frac{\prod_{i=1}^k |C_i|}{|G|}\sum_{\c\in \Irr(G)} \frac{\c(c_1)\cdots \c(c_k)\c(g^{-1})}{\c(1)^{k-1}}
\]
expressing the number of ways of writing an element $g \in G$ as a product $x_1x_2\cdots x_k$ of elements $x_i\in C_i$, where $C_i = c_i^G$ are $G$-conjugacy classes of elements $c_i$, $1 \leq i \leq k$, and the sum is over the set $\Irr(G)$ of all irreducible characters of $G$ (see \cite[10.1]{AH}).  This connection is sometimes a starting point for such applications;
it has been particularly exploited for almost simple (or quasisimple) groups $G$.

Another classical formula involving character ratios goes back to Frobenius in 1896 \cite{F}.
It asserts that, for any finite group $G$, the number $N(g)$ of ways to express an element $g \in G$
as a commutator $[x,y]$ ($x,y \in G$) satisfies
\[
N(g)= |G| \sum_{\chi \in \Irr(G)} \frac{\chi(g)}{\chi(1)}.
\]
This formula is widely used, and served (together with character bounds) as an important tool
in the proof of Ore's conjecture \cite{LOST}.

We are particularly interested in so called {\it exponential character bounds}, namely bounds of the
form
\[
|\chi(g)| \le \chi(1)^{\a_g},
\]
sometimes with a multiplicative constant, holding for {\it all} characters $\chi \in \Irr(G)$, where $0 \le \a_g \le 1$ depends on the group element $g \in G$. Obviously, if $g$ is central in $G$, then we must have $\a_g =1$, but for most elements $g$ we aim to find
$\a_g < 1$ which is as small as possible. One advantage of exponential character bounds is that they imply the inequality $|\frac{\c(g)}{\c(1)}| \le \chi(1)^{-(1-\a_g)}$, so the upper bound on the character ratio becomes smaller
as the character degree grows.

The first exponential character bound was established in 1995 for symmetric groups $\SSS_n$ by Fomin and Lulov
\cite{FL}.  They show that, for permutations $g \in \SSS_n$ which are products of $n/m$
cycles of length $m$ and for all characters $\chi \in \Irr(\SSS_n)$ we have
\begin{equation}\label{fomlul}
|\chi(g)| \le c(m) n^{\frac{1}{2}-\frac{1}{2m}} \chi(1)^{\frac{1}{m}},
\end{equation}
for a suitable function $c: \N \to \N$.

In \cite{fuchs1} this bound and some extensions of it were applied in various contexts, including
the theory of Fuchsian groups.
Subsequently, exponential character bounds which hold for {\it all} permutations $g \in \SSS_n$ and which are
essentially best possible were established in 2008 in \cite{LarS}, with applications to a range of problems: mixing times of random walks, covering by powers of conjugacy classes, as well as probabilistic and combinatorial properties of word maps.

Can we find good exponential character bounds for groups of Lie type? This problem has turned out to be
quite formidable; it has been considered by various researchers over the past two decades, and various approaches have been attempted, but it is only in this paper that strong (essentially best possible) such bounds are established.

The first significant bound on character ratios for groups of Lie type was obtained in 1993 by Gluck \cite{Gl1},
who showed that $\frac{|\c(g)|}{\c(1)} \le Cq^{-1/2}$ for any non-central element $g \in G(q)$, a group of Lie type over $\F_q$, and any non-linear irreducible character $\c$ of $G(q)$, where $C$ is an absolute constant. In \cite{gluck}, he proved a bound of the form
\[
\frac{|\c(g)|}{\c(1)} \le \c(1)^{-\g/n},
\]
when $G(q)$ is a classical group with natural module $V = \F_q^n$ of dimension $n$, and $\g = \g(q,d)$ is a positive
real number depending on $q$ and on $d = \dim [V,g]$, the dimension of the commutator space of $g$ on $V$.
While this result provides an exponential character bound $|\c(g)| \le \c(1)^{\a_g}$,
the exponent $\a_g = 1 -\g/n$ is not explicit, and in the general case we have $\g(q,d) \le 0.001$, so
$\a_g \ge 1- \frac{1}{1000n}$, which is very close to $1$.

An explicit character bound for finite classical groups, with natural module $V = \F_q^n$,
in terms of the {\it support} $\supp(g)$ of the element $g$ was obtained in \cite[4.3.6]{LST}: namely,
\begin{equation}\label{support}
\frac{|\c(g)|}{\c(1)} < q^{-\sqrt{\supp(g)}/481},
\end{equation}
where $\supp(g)$ is the codimension of the largest eigenspace of $g$ on $V \otimes_{\F_q}\overline\F_q$.
These results have applications to covering, mixing times and word maps.

In this paper we obtain asymptotically much stronger bounds for character ratios of finite groups of Lie type in good characteristic (this restriction comes
from the fact that our proof relies on certain results in the Deligne-Lusztig theory, which currently are only known to hold in good characteristics). In fact we provide the first explicit exponential character bounds for groups of
Lie type, and show that these bounds are asymptotically optimal in many cases.

These character bounds lead to several new results on random walks and covering by products of conjugacy classes
that are far stronger than previously known such results. Further applications to the theory of representation varieties of Fuchsian groups and probabilistic generation of groups of Lie type will be given in a sequel to this paper \cite{lst2}.

We also prove the first bounds on character ratios for Brauer characters, for the groups $SL_n(q)$ and $GL_n(q)$,
and in characteristics coprime to $q$.

We now describe our results. Throughout the paper, let $\KK$ be an algebraically closed field of characteristic $p$,
$\GC$ a connected reductive algebraic group over $\KK$, $F:\GC \to \GC$ a Frobenius endomorphism,
and $G = \GC^F$.
For a subgroup $X$ of $\GC$ write $X_{{\rm {unip}}}$ for the set of {\it non-identity} unipotent elements
of $X$. For a fixed $F$, a Levi subgroup $\LC$ of $\GC$ will be called {\it split}, if it is an $F$-stable
Levi subgroup of an $F$-stable proper parabolic subgroup of $\GC$. For an $F$-stable Levi subgroup
$\LC$ of $\GC$ and $L=\LC^F$, we define
$$\a(L) := \hbox{ {\rm {max}} }_{u \in L_{{\rm {unip}}}} \frac{\dim u^\LC}{\dim u^\GC},~~
    \a(\LC) := \hbox{ {\rm {max}} }_{u \in \LC_{{\rm {unip}}}} \frac{\dim u^\LC}{\dim u^\GC}$$
if $\LC$ is not a torus, and $\a(L)=\a(\LC) := 0$ otherwise.

\begin{thm}\label{main1}
There exists a function $f:\N \to \N$ such that the following statement holds.
Let $\GC$ be a connected reductive algebraic group such that $[\GC,\GC]$ is simple of rank $r$ over a field of
good characteristic $p > 0$. Let $G := \GCF$ for a Frobenius endomorphism $F: \GC \to \GC$.
Let $g \in G$ be any element
such that $\CB_G(g) \leq L := \LC^F$, where $\LC$ is a split Levi subgroup of $\GC$. Then,
for any character $\chi \in \Irr(G)$ and $\al:=\al(L)$, we have
$$|\c(g)| \leq f(r)\c(1)^\a.$$

\end{thm}

\begin{rem}\label{rems}
{\em (i) The $\al$-bound in Theorem \ref{main1} is sharp  in several cases -- see Example \ref{sharpness}. In fact,
this $\al$-bound is {\it always} sharp in the case of $GL_n(q)$ and
$SL_n(q)$, by Theorem \ref{main2}.

(ii) If $r \geq 9$ and $q \geq r^2+1$, then the function $f(r)$ in Theorem \ref{main1} can be chosen to be
$2^{2r+\sqrt{2r}+3} \cdot (r!)^2$ (with the main term being the square of the largest order of the
Weyl group of a simple algebraic group of rank $r$) -- see Proposition \ref{frbd}. Moreover,
$\al \lesssim_r 1-1/r$ by Theorem \ref{ratio} and $\chi(1) \geq q^r/3$ if $\chi(1) > 1$ by \cite{LSe}, hence
Theorem \ref{main1} yields $|\chi(g)| \lesssim_r  \chi(1)^{\a+1/2r} \lesssim_r \chi(1)^{1-1/2r}$ if $q > r^{4r}$; in fact,
$\chi(1) \geq q^{r^2/2}$ for most of $\chi \in \Irr(G)$, for which the bound becomes
$|\chi(g)| \lesssim_r \chi(1)^{\a+1/3r} \lesssim_r \chi(1)^{1-2/3r}$ if $q > r^{12}$. (Here, we say
that $f_1(x) \lesssim_x f_2(x)$ for two functions
$f_1,f_2 : \R \to \R_{\geq 0}$ if $\limsup_{x \to \infty} f_1(x)/f_2(x) \leq 1$.)

(iii) Although the aforementioned choice of $f(r)$ in Theorem \ref{main1} can be improved, Example \ref{sharpness}(vi) shows that
$f(r)$ should be at least the largest degree of complex irreducible characters of the Weyl group
$W(\GC)$ of $\GC$, which can be quite close to $|W(\GC)|^{1/2}$.
In particular, choosing $\GC$ of type $A_r$ and applying \cite{LoS}, \cite{VK}, we get
$f(r) > e^{-1.283\sqrt{r+1}}\sqrt{(r+1)!}$.}
\end{rem}

Note that Theorem \ref{main1} and its various consequences also apply for finite twisted groups
of Lie type.

\begin{thm}\label{main2}
In the notation of Theorem \ref{main1}, there is a constant $C_n >0$ depending only on
$n$ such that the following statement holds. For $G = \GC^F = GL_n(q)$ with
$q \geq C_n$ and for any split Levi subgroup $\LC$ of $\GC$, there is a semisimple element
$g \in G$ and a unipotent character $\chi \in \Irr(G)$ such that $\CB_G(g) = L = \LC^F$ and
$$\chi(g) \geq \frac{1}{4} \chi(1)^{\a(L)}.$$
The same conclusion holds for $SL_n(q)$, if for instance we choose $q$ so that
$q-1$ is also divisible by
$(n!)^n$.
\end{thm}

\medskip
In the case of $GL_n(q)$ and $SL_n(q)$
we can also prove a version of Theorem \ref{main1} for Brauer characters in cross-characteristic.

\begin{thm}\label{main1b}
There exists a function $h:\N \to \N$ such that the following statement holds.
Let $\GC = GL_n$ or $SL_n$ be an algebraic group over a field of characteristic
$p > 0$ and  $F: \GC \to \GC$ a Frobenius endomorphism,
such that $G = \GC^F \cong GL_n(q)$ or $SL_n(q)$.
Let $\ell = 0$ or a prime not dividing $q$.
Let $g \in G$ be any $\ell'$-element
such that $\CB_G(g) \leq L := \LC^F$, where $\LC$ is a split Levi subgroup of $\GC$.
Then for any irreducible $\ell$-Brauer character $\varphi$ of $G$ and $\a := \a(L)$,
we have
$$|\varphi(g)| \leq h(n)\varphi(1)^\a.$$
\end{thm}

The above results do not cover, for instance, the case where $g \in \GCF$ is a unipotent element.
However, we have been able to obtain a complete result covering all elements in $GL_n(q)$ and $SL_n(q)$:

\begin{thm}\label{main1c}
There is a function $h:\N \to \N$ such that the following statement holds. For any $n \geq 5$, any prime
power $q$, any irreducible complex character $\chi$ of $H:= GL_n(q)$ or $SL_n(q)$, and any
non-central element $g \in H$,
$$|\chi(g)| \leq h(n) \cdot \chi(1)^{1-\frac{1}{2n}}.$$
\end{thm}

For the remaining groups of Lie type, character bounds, which work for arbitrary elements $g \in \GCF$, and are
weaker than the one in Theorem \ref{main1} but
asymptotically stronger than the ones in \cite{gluck} and \cite{LST}, will be proved in a sequel to this paper.

\smallskip

To be able to apply Theorem \ref{main1} we need information on the values of $\a(L) \leq \a(\LC)$.
For classical groups, we prove the following upper bound.

\begin{thm} \label{ratio} If $\GC$ is a classical algebraic group over $\KK$ in good characteristic, and $\LC$ is a Levi subgroup of $\GC$, then
\[
\a(\LC) \le {1\over 2}\left( 1+ {{\dim \LC} \over {\dim \GC}}\right).
\]
\end{thm}

For exceptional types we obtain fairly complete information.

\begin{thm}\label{alphaexcep}
If $\GC$ is an exceptional algebraic group in good characteristic, the values of $\a(\LC)$ for (proper, non-toral) Levi subgroups $\LC$ are as in Table $\ref{extab}$.
\end{thm}

In Table \ref{extab}, for $\GC = F_4$ or $G_2$ the symbols $\tilde A_1, \tilde A_2$ refer to Levi subsystems consisting of short roots. For $\GC = E_7$, there are two Levi subgroups $A_5$ and $A_5'$: using the notation for the fundamental roots $\a_i\,(1\le i\le 7)$ as in \cite{bour}, these are the Levi subgroups with fundamental roots $\{\a_i : i=1,3,4,5,6\}$ and $\{\a_i : i=2,4,5,6,7\}$ respectively. The notation  $\triangleright A_4$, for instance, means that
$\LC'=[\LC,\LC]$ has a simple factor of type $A_4$.

\begin{table}[h!]
\caption{$\a$-values for exceptional groups} \label{extab}
\vspace{-7mm}
\[
\begin{array}{r|cccccccccc}
\hline
&&&&&&&&&& \\
\vspace{-8mm}\\
\GC=E_8,\,\LC '= & E_7 & D_7 & \LC'\triangleright E_6 & D_6 & A_7 & \triangleright D_5 & \triangleright A_6 & \triangleright A_5 &
\triangleright D_4 & \hbox{rest} \\
\a(\LC)= & \frac{17}{29} & \frac{9}{23} & \frac{11}{29} & \frac{9}{29} & \frac{15}{56} &  \frac{7}{29} & \frac{5}{23}
&  \frac{4}{23} &  \frac{5}{29} &  \le \frac{1}{6}  \\
\vspace{-3mm}\\
\hline
&&&&&&&&&& \\
\vspace{-8mm}\\
\GC=E_7,\,\LC '= & E_6 & D_6 & \LC'\triangleright D_5 & A_6 & A_5 & \triangleright A_5' & \triangleright D_4 & \triangleright A_4 &
\triangleright A_3 & \hbox{rest} \\
\a(\LC)= & \frac{11}{17} & \frac{5}{9} & \frac{7}{17} & \frac{5}{13} & \frac{4}{13} &  \frac{1}{3} & \frac{5}{17}
& \le \frac{1}{4} &  \le \frac{1}{5} &  \le \frac{1}{6}  \\
\vspace{-3mm}\\
\hline
&&&&&&&&&& \\
\vspace{-8mm}\\
\GC=E_6,\,\LC '= & D_5 & A_5 & D_4 & \LC'\triangleright A_4 & \triangleright A_3 & \triangleright A_2 & A_1^k &&&\\
\a(\LC)= & \frac{7}{11} & \frac{1}{2} & \frac{5}{11} & \frac{3}{8} & \frac{3}{11} & \le \frac{7}{27} & \le \frac{3}{20}&&& \\
\vspace{-3mm}\\
\hline
&&&&&&&&&& \\
\vspace{-8mm}\\
\GC=F_4,\,\LC ' =& B_3 & C_3 & A_2\tilde A_1,A_2 & \tilde A_2A_1 & \tilde A_2 & A_1\tilde A_1 & A_1 & \tilde A_1 &&\\
\a(\LC)= & \frac{1}{2} & \frac{7}{15} & \frac{1}{4} & \frac{2}{9} & \frac{1}{5} & \frac{1}{7} & \frac{1}{8} & \frac{1}{11}&& \\
\vspace{-3mm}\\
\hline
&&&&&&&&&& \\
\vspace{-8mm}\\
\GC=G_2,\,\LC '= & A_1 & \tilde A_1&&&&&&&& \\
\a(\LC)= & \frac{1}{3} & \frac{1}{4}&&&&&&&& \\
\vspace{-3mm}\\
\hline
\end{array}
\]
\end{table}

We can now easily deduce the following.

\begin{cor}\label{keycor}
Let $\GC$, $G = \GC^F$, and $f$ be as in Theorem $\ref{main1}$. Suppose $y\in G$ is a (semisimple) element
such that $\CB_G(y) = \LC^F$, where $\LC$ is a
split Levi subgroup of $\GC$. Then for any non-linear $\chi \in {\rm Irr}(G)$,
\[
|\c(y)| \le f(r) \, \c(1)^{1-\frac{1}{2}\frac{\dim y^{\GC}}{\dim \GC}}.
\]
\end{cor}

Next, we establish a new strong bound on character ratios given the support (which is defined right after \eqref{support})
of the semisimple part of the ambient element.


\begin{thm}\label{supp}
Assume $\GC = SL_n(\overline\F_q)$ with $n \geq 2$, $Sp_n(\overline\F_q)$ with $n \geq 4$, or $Spin_n(\overline\F_q)$
with $n \geq 7$, all in good characteristic, and define
$$c := c(\GC) = \left\{ \begin{array}{ll} (r+1)/(2r+4), & \GC = SL_{r+1},\\
    r/(4r+2), & \GC = Sp_{2r},\\
    r/(4r-2), & \GC = Spin_{2r},\\
    1/4, & \GC = Spin_{2r+1}. \end{array} \right.$$
Let $G = \GC^F = G(q)$ and $f$ be as in Theorem \ref{main1}, and let $g \in G$ be any element
such that its semisimple part $y$ has centralizer $\CB_G(y) = \LC^F$,
where $\LC$ is a split Levi subgroup of $\GC$.
Then, for any non-linear $\chi \in {\rm Irr}(G)$,
\[
\frac{|\c(g)|}{\c(1)} \leq 3f(r) \, q^{-c\cdot\supp(y)}.
\]
\end{thm}

In particular it follows that $\frac{|\c(y)|}{\c(1)} \leq 3f(r) \, \c(1)^{-c\cdot\supp(y)}$, and that
for any $\e > 0$, $r \ge r(\e)$ and $q$ larger than a suitable function of $r$, we have
\[
\frac{|\c(y)|}{\c(1)} \le q^{-(b-\e)\cdot\supp(y)},
\]
where $b = 1/2$ in the $SL_{r+1}$ case and $b = 1/4$ in the other cases.

Theorem \ref{supp} and its consequences considerably improve the bound (\ref{support}) from \cite[4.3.6]{LST}
for elements as above.

\medskip

We also obtain more precise character bounds for $GL_n$. To state them we need some notation.
For positive integers $n_1, \ldots , n_m$ define
\[
\beta(n_1, \ldots , n_m) :=  \max \frac{\sum_{i=1}^m (n_i^2 - \sum_{j=1}^n a_{ij}^2)}
{n^2 - \sum_{j=1}^n (\sum_{i=1}^m a_{ij})^2},
\]
where $n = n_1 + \ldots + n_m$ and the maximum is taken over all non-negative integers $a_{ij}$ ($1\le i \le m$,
$1 \le j \le n$) satisfying
$$\sum_{j=1}^n a_{ij} = n_i, ~~a_{i1} \geq a_{i2} \geq \ldots \geq a_{in}, ~~1 \leq i \leq m,~~~\max_{1 \leq i \leq m}a_{i2} > 0,$$
if $\max_{1 \leq i \leq m}n_i \geq 2$, and let $\beta(1,1, \ldots,1) = 0$.


\begin{thm} \label{GL} Let $G = GL_n(q)$
and let $L \le G$ be a Levi subgroup of the form $L = GL_{n_1}(q) \times \cdots \times GL_{n_m}(q)$,
where $n_i \ge 1$ and $\sum_{i=1}^m n_i =n$. Let $n_{i_0} = \max_{1 \leq i \leq m}n_i$. Then
\[
  \frac{n_{i_0}-1}{n-t} \leq \a(L) = \beta(n_1, \ldots, n_m) \leq \frac{n_{i_0}}{n}
\]
if $n_{i_0} \geq 2$ and $t$ is the number of $1 \leq j \leq m$ such that $n_j = n_{i_0}$,
and $\a(L) = \beta(n_1, \ldots ,n_m) = 0$ if $n_{i_0} = 1$. Consequently, for every $g \in G$ with $\CB_G(g) \le L$ and every
$\chi \in \Irr(G)$ we have
\[
|\chi(g)| \le f(n-1) \chi(1)^{\beta(n_1, \ldots , n_m)},
\]
where $f: \N \to \N$ is the function specified in Theorem \ref{main1}.
\end{thm}


Suppose now that $m$ divides $n$ and $n_1 = \ldots = n_m = n/m >1$.
Then we can show that $\beta(n_1, \ldots , n_m) = \frac{1}{m}$, so we immediately obtain the following.

\begin{cor} \label{FL} Let $G = GL_n(q)$ where $q$ is a prime power. Let $m < n$
be a divisor of $n$ and let $L \le G$ be a Levi subgroup of the form $L = GL_{n/m}(q)^m$.
Let $g \in G$ with $\CB_G(g) \le L$. Then we have
\[
|\chi(g)| \le f(n-1) \chi(1)^{\frac{1}{m}}
\]
for all characters $\chi \in \Irr(G)$,
where $f: \N \to \N$ is the function specified in Theorem \ref{main1}.
\end{cor}

Example \ref{sharpness} again shows that the exponent $1/m$ in Corollary \ref{FL} is sharp. In general, Theorem \ref{GL} determines
$\a(L)$ up to within $1/n$. It is reasonable to conjecture
that, under the hypotheses of Theorem \ref{GL}, $\a(L) = (n_{i_0}-1)/(n-t)$.
This conjecture is confirmed in Theorem \ref{GL2} for the case $m=2$ (as well as in the cases, where either $n \leq 8$,
or $m \leq 4$ and $n \leq 13$, by direct calculation).

The bound in Theorem \ref{GL}  and some variations on it have applications to Fuchsian groups
(see \cite{lst2}). Corollary \ref{FL} may be regarded as a Lie analogue of the Fomin-Lulov
character bound (\ref{fomlul}) for $\SSS_n$ mentioned before.

\vspace{4mm}
We now present some applications of the above results to the theory of
{\it mixing times} for random walks on finite quasisimple groups of Lie type corresponding to conjugacy classes.
Let $G = G(q)$ be such a group, let $y\in G$ be a non-central element, and let $C = y^G$, the conjugacy class of $y$.
Consider the random walk on the corresponding Cayley graph starting at the identity, and at each step moving from a vertex $g$ to a neighbour $gs$, where $s \in y^G$ is chosen uniformly at random. Let $P^t(g)$ be the probability of reaching the vertex $g$ after $t$ steps. The mixing time of this random walk is defined to be the smallest integer $t = T(G,y)$ such that $||P^t-U||_1< \frac{1}{e}$, where $U$ is the uniform distribution and $||f||_1 = \sum_{g\in G} |f(g)|$ is the $l_1$-norm.

Mixing times of such random walks have been extensively studied since the pioneering work of Diaconis and
Shashahani \cite{DS} on the case $G = \SSS_n$ and $C$ the class of transpositions in $\SSS_n$.
Additional results on random walks in symmetric and alternating groups have been obtained in various
papers, see for instance \cite{Roi}, \cite{V}, \cite{LP} and \cite{LarS}. The latter paper obtains essentially
optimal results on mixing times in these groups.

However, if we turn from symmetric groups to finite groups $G$ of Lie type, good estimates on mixing
times have been obtained only in very few cases.  Hildebrand \cite{Hil} showed that the mixing time for
the class of tranvections in $SL_n(q)$ is of the order of $n$.  In \cite{chardeg} it is shown that
if $y \in G$ is a regular element, then the mixing time $T(G,y)$ is $2$ when $G \ne PSL_2(q)$ is large.
In \cite{S1} it is proved that, if $G$ is any finite simple group, then for a random $y \in G$ we have
$T(G, y) = 2$ (namely, the latter equality holds with probability tending to $1$ as $|G| \to \infty$).
Other than that, the mixing times $T(G,y)$ for groups $G$ of Lie type remain a mystery.

The next result contains bounds for mixing times, and also (in parts (I)(a) and (II)) for the number of steps required so that $P^t$ is close to $U$ in the $l_\infty$-norm, which is stronger than the $l_1$-norm condition for mixing time (and also implies that the random walks hits all elements of $G$). Here we define
$||f||_\infty = |G|\,{\rm max}_{x\in G}|f(x)|$, and say that $C^t = G$ {\it almost uniformly pointwise} as $q\rightarrow \infty$ if  $||P^t-U||_\infty \rightarrow 0$ as $q\rightarrow \infty$.

We denote by $h:=h(\GC)$, the Coxeter number of $\GC$, defined by
\[
h(\GC) = \frac{\dim{\GC}}{r} -1,
\]
where $r$ is the rank of $\GC$. Note that $h \ge 2$ and that $h \to \infty$ as $r \to \infty$.

\begin{thm}\label{mix}
Suppose $\GC$ is a  simple algebraic group in good characteristic, and $G = G(q) = \GC^F$ is a finite quasisimple group over $\F_q$.
Let $y \in G$ be such that $\CB_G(y) \le  L$, where $L = \LC^F$ for a split Levi subgroup $\LC$ of $\GC$. Write $C = y^G$.
\begin{itemize}
\item[{\rm (I)}] Suppose $\GC$ is of classical type.
\begin{itemize}
\item[{\rm (a)}] If $t > (4+\frac{4}{h})\frac{\dim \GC}{\dim \GC- \dim \LC}$, then
$C^t = G$ almost uniformly pointwise as $q \rightarrow \infty$.
In particular, $C^t = G$ for sufficiently large $q$.
\item[{\rm (b)}] The mixing time $T(G,y) \le \lceil (2+\frac{2}{h})\frac{\dim \GC}{\dim \GC - \dim \LC}\rceil$ for large $q$.
\end{itemize}
\item[{\rm (II)}] Suppose $\GC$ is of exceptional type. Then $C^6 = G$ almost uniformly pointwise as $q \rightarrow \infty$, and the mixing time $T(G,y) \le 3$.
\end{itemize}
\end{thm}

\no
{\bf Remarks } (i) Note that the multiplicative constants above are very small.
For example, $2 + \frac{2}{h} \le 3$ and it tends to $2$ as $r \to \infty$.

(ii) The constant $2+\frac{2}{h}$ in part
I(b) of Theorem \ref{mix}  is best possible for some classes, for example homologies $y = \hbox{diag}(\mu I_{n-1},\la)$ in $G = SL_n(q)$
(where $\mu,\la \in \F_q^\times$ and $\mu \neq \la$ and so $q \geq 3$), for which the bound given by part (I)(b) is $T(G,y) \le n+3$ and
for which the mixing time is at least $n$ by Lemma \ref{mix-subset}(ii).

(iii) The bound  $T(G,y) \le 3$ for exceptional groups in (II) is best possible for many classes -- namely, those classes for which $\dim y^\GC$ is smaller than $\frac{1}{2}\dim \GC$. For such classes, $|y^G|^2 < |G|/2$
for large $q$, so the mixing time cannot be 2 by Lemma \ref{mix-subset}(i).

\vspace{4mm}




Theorem \ref{mix}(I)(b) implies the following linear bounds for classical groups.

\begin{cor}\label{linearbd} Let $G = \GC^F$ be a quasisimple classical group over $\F_q$, where $\GC$ is simple of rank $r$ over $\F_q$,
and let $y \in G$ be as in Theorem $\ref{mix}$. Then for large $q$,
\begin{enumerate}[\rm(i)]
\item the diameter ${\rm diam}(G,y^G) \le 2r+4$, and
\item the mixing time $T(G,y) \le r+2$.
\end{enumerate}
\end{cor}

A linear bound for the diameter (of the order of $40r$), which holds for all non-central conjugacy classes, can be found in \cite{LLi}.

Using Theorem \ref{main1c} we can obtain such a bound for {\it all} conjugacy classes in $SL_n(q)$:

\begin{cor}\label{linearbdsl} Let $G= SL_n(q)$,  let $x$ be an arbitrary non-central element of $G$ and let $C = x^G$.
\begin{itemize}
\item[{\rm (i)}] If $t>4n+4$, then $C^t = G$ almost uniformly pointwise as $q\rightarrow \infty$.
\item[{\rm (ii)}] The mixing time $T(G,x) \le 2n+3$ for large $q$.
\end{itemize}
\end{cor}

Note that \cite[Theorem 1]{lev} shows that  $C^n=G$ for any nontrivial conjugacy class $C$ of $G=PSL_n(q)$, where $n\ge 3, q\ge 4$.

We can also use Theorem \ref{gl-uni} (or rather its corollary \ref{sl-uni}) to obtain a better bound for unipotent elements of $SL_n(q)$.

\begin{thm}\label{mix-gluni}
Let $G = SL_n(q)$ and let $u$ be a non-identity unipotent element in $G$. Write $C = u^G$.
\begin{itemize}
\item[{\rm (i)}] If $t > 2n$, then
$C^t = G$ almost uniformly pointwise as $q \rightarrow \infty$.
In particular, $C^t = G$ for sufficiently large $q$.

\item[{\rm (ii)}] The mixing time $T(G,u) \le n$ for sufficiently large $q$.
\end{itemize}
\end{thm}


One can compare part (ii) of the above theorem with Hildebrand's result \cite{Hil} for transvections, where he proves that for $n$ varying, the mixing time for the class of tranvections in $SL_n(q)$ is of the order of $n$. In our case
$n$ may still vary, but $q$ should be much larger than $n$.  The coincidence of values seems striking.


It is interesting to compare the mixing time $T(G,y)$ with the {\it covering number} $cn(G,C)$
of the conjugacy class $C = y^G$, defined as the minimal $t$ for which $C^t = G$.
It is known that there is an absolute constant $b$ such that for any conjugacy class
$C \ne \{ 1 \}$ of any finite simple group $G$ we have
\[
\frac{\log |G|}{\log |C|} \le cn(G,C) \le b \frac{\log |G|}{\log |C|}.
\]
Indeed the first inequality is trivial, while the second is \cite[1.2]{lishdiam}.

It is easy to see that, with the above notation,
\begin{equation}\label{time}
\frac{\log |G| + \log (1-e^{-1})}{\log |C|} \le T(G,y).
\end{equation}
Indeed, this follows from Lemma \ref{mix-subset}.

It is conjectured in \cite[4.3]{S3} that there is an absolute constant $c$ such that
for any finite simple group $G$ of Lie type and any non-identity element $y \in G$ we have
\begin{equation}\label{1.1}
T(G, y) \le c \frac{\log |G|}{\log |C|},
\end{equation}
where $C = y^G$.

Note that this statement does not hold for alternating groups $G$ (take $y \in G$ to be
a cycle of length around $n/2$ -- then $\frac{\log |G|}{\log |C|}$ is bounded, while
$T(G,y)$ is of the order of $\log n$).

The above conjecture is related to an older conjecture posed by Lubotzky in \cite[p.179]{lub}.
Lubotzky conjectured that, if $G$ is a finite simple group and $C$ is a non-trivial conjugacy class
of $G$, then the mixing time of the Cayley graph $\Gamma(G, C)$ of $G$ with $C$ as a generating set is
linearly bounded above in terms of the diameter of $\Gamma(G, C)$. Since this diameter is exactly
the covering number $cn(G,C)$, this conjecture (combined with the more recent upper bound on
$cn(G,C)$ mentioned above) implies conjecture (\ref{1.1}).

Applying Theorem \ref{mix} we are able to prove the above conjectures in many interesting cases.
Let $\GC$ be a simple algebraic group in good characteristic, and $G = G(q) = \GC^F$ a finite quasisimple group
over $\F_q$. We say that a non-central element $y \in G$ is {\it nice} if $\CB_G(y) = L$, where $L = \LC^F$ for a split Levi subgroup $\LC$ of $\GC$. Note that (non-central) split semisimple elements are nice.

\begin{cor}\label{lub}
Let $\GC$, $G(q)$ be as above, and suppose $q$ is large (given $\GC$).
Then Conjecture (\ref{1.1}) holds for all nice elements $y$ of the quasisimple group $G(q)$.
In particular, the conjecture holds for all split semisimple elements of $G(q)$.
\end{cor}

Indeed, this readily follows from part I(b) of Theorem \ref{mix}, with a very small constant
$c$ (around $3$).

Conjecture (\ref{1.1}) and Corollary \ref{lub} suggest a distinctive difference between mixing times for
$\SSS_n$ as opposed to classical groups $Cl_n(q)$.


Our final result essentially determines the mixing time $T(G,y)$ in terms of the support of $y$ as follows (recall
the notation $f_1(x) \lesssim_x f_2(x)$ from Remark \ref{rems}).

\begin{thm}\label{mix-supp}
Assume $\GC = SL_n(\overline\F_q)$ with $n \geq 2$, $Sp_n(\overline\F_q)$ with $n \geq 4$, or $Spin_n(\overline\F_q)$
with $n \geq 7$, and define
$$r' := r'(\GC) = \left\{ \begin{array}{ll} r(2r+4)/(r+1), & \GC = SL_{r+1},\\
    4r+2, & \GC = Sp_{2r},\\
    4r-2, & \GC = Spin_{2r},\\
    4r, & \GC = Spin_{2r+1}. \end{array} \right.$$
Let $G = \GC^F = G(q)$ and $f$ be as in Theorem \ref{main1}, and let $g \in G$ be any element
such that its semisimple part $y$ has centralizer $\CB_G(y) = \LC^F$,
where $\LC$ is a split Levi subgroup of $\GC$.
Suppose $q$ is large enough (given $r$).
Then we have
\[
T(G,g) \le \lceil (2+\frac{2}{h})r'/\supp(y)\rceil.
\]
Furthermore, we have
\[
\frac{1}{2} r'/\supp(y) \lesssim_{|G|} T(G,y) \le \lceil (2+\frac{2}{h})r'/\supp(y)\rceil.
\]
\end{thm}

Thus, under the above conditions, the mixing time $T(G,y)$ is essentially $n/\supp(y)$
(up to a small multiplicative constant).

\section{Character bounds: Proof of Theorem \ref{main1}}\label{sec:main1}
Throughout this section, let $\GC$ be a connected reductive algebraic group over a field of characteristic
$p > 0$, $F:\GC \to \GC$ a Frobenius endomorphism, and let $G := \GC^F$. We will say
$\GCF$ is defined over $\F_q$, if $q$ is the common absolute value of eigenvalues of $F$ acting
on $X(\TC)\otimes \R$, where $X(\TC)$ is the character group of an $F$-stable (maximally split) maximal torus
$\TC$ of $\GC$.

First we prove the following statement concerning Harish-Chandra restriction.

\begin{prop}\label{fixed}
Suppose that $g  \in G$ is such that $\CB_G(g) \leq \LC^F$, where $\LC$ is an $F$-stable Levi
subgroup of an $F$-stable parabolic subgroup
$\PC = \UC\LC$ of $\GC$ with unipotent radical $\UC$.
Let $\ell = 0$ or a prime not dividing $p|g|$,
$\F$ an algebraically closed field of characteristic $\ell$, and let
$\varphi$ be the Brauer character of some $\F G$-module $V$. Also, let $\psi$ denote the Brauer
character of the $\LC^F$-module $\CB_V(\UC^F)$.
Then
$$\varphi(g) = \psi(g).$$
\end{prop}

\pf
(a) Write $L := \LC^F$, $P := \PC^F$, and $U := \UC^F$.
First we handle the case $\ell = 0$. Consider the map $f:U \to U$ given by $f(u) = g^{-1}ugu^{-1}$. Then,
for $u, v \in U$ we have that
$$f(u) = f(v) \Leftrightarrow v^{-1}u \in U \cap \CB_G(g) \subseteq U \cap L = 1 \Leftrightarrow u = v.$$
Thus the map $f$ is injective, and so bijective. Hence, when $u$ runs over $U$, $ugu^{-1}$ runs over
the elements of $gU$, each element once:
$$\{ugu^{-1} \mid u \in U\} = gU.$$
 Now we decompose $V = \CB_V(U) \oplus [V,U]$ as a $P$-module
(note that $P = \NB_G(U)$), and let $\Phi= \diag(\Phi_1,\Phi_2)$ denote the representation of $P$
with respect to some basis respecting this decomposition. In particular, no irreducible constituent of $(\Phi_2)|_U$ is
trivial, and so $\sum_{u \in U}\Phi_2(u) = 0$. It follows that
$$\sum_{u \in U}\Phi(ugu^{-1}) = \sum_{u \in U}\Phi(gu) = \Phi(g)\sum_{u \in U}\Phi(u) = $$
$$   = \diag(\Phi_1(g)\sum_{u \in U}\Phi_1(u),\Phi_2(g)\sum_{u \in U}\Phi_2(u))
       = \diag(|U|\Phi_1(g),0).$$
Taking the trace of both sides, we obtain $|U|\varphi(g) = |U|\psi(g)$, as stated.


\smallskip
(b) For the modular case $\ell > 0$, let $\chi^\circ$ denote the restriction of any complex character $\chi$ of
$G$ or $P$ to $\ell'$-elements. It is well known, see e.g \cite[Theorem 15.14]{Is}, that any Brauer character of $G$ is
a $\Z$-combination of $\chi^\circ$ with $\chi \in \Irr(G)$. It follows that (in the Grothendieck group of $\F G$-modules)
we can write $V = V_1-V_2$, where $V_1$ and $V_2$ are some reductions modulo $\ell$ of $\C G$-modules
$W_1$ and $W_2$ affording complex characters $\chi_1$ and $\chi_2$. Since $\ell \neq p$,
$\CB_V(U) = \CB_{V_1}(U)-\CB_{V_2}(U)$ in the Grothendieck group of $\F P$-modules. Now $g \in P$,
$\varphi(g) = \chi_1(g)-\chi_2(g)$, and the statement follows
by applying the results of (a) to $W_1$ and $W_2$.
\hal

Recall that the complex irreducible characters of $G = \GC^F$ can be partitioned into {\it Harish-Chandra series}, see
\cite[Chapter 9]{C}. We refer to \cite{C} and \cite{DM} for basic facts on {\it Harish-Chandra restriction}
$\SR^G_L$ and {\it Harish-Chandra induction} $R^G_L$.

\begin{prop}\label{hc}
There is a constant $A = A(r)$ depending only on the semisimple rank $r$ of $\GC$ with the following property.
Suppose that $\chi \in \Irr(G)$ is such that $\SR^G_L(\chi) \neq 0$ for $L = \LC^F$, where $\LC$ is a split Levi
subgroup of $\GC$. Then the total number of irreducible constituents of
the $L$-character $\SR^G_L(\chi)$ (with counting multiplicities) is at most $A$. In fact, if $[\GC,\GC]$ is simple then
one can choose $A = W(r)^2$, where $W(r)$ denotes the largest order of the Weyl group of  a simple algebraic group of rank $r$.
\end{prop}

\pf
Since $\SR^G_L(\chi) \neq 0$, $\chi$ is not cuspidal. By \cite[Proposition 9.3.1]{C}, we may assume that
$\LC$ is a standard $F$-stable Levi subgroup of a standard $F$-stable parabolic subgroup $\PC=\UC\LC$  of $\GC$.
Suppose that $\chi$ belongs to the Harish-Chandra series
labeled by a standard Levi subgroup $L_1$ and a cuspidal character $\psi \in \Irr(L_1)$. Here,
$L_1 = \LC_1^F$, where $\LC_1$ is a split Levi subgroup of $\GC$,
and $\chi$ is an irreducible constituent of $R^G_{L_1}(\psi)$.

Suppose now that $\eta$ is any irreducible constituent of $\SR^G_L(\chi)$, and let $\eta$ belongs to
the Harish-Chandra series labeled by a standard Levi subgroup $L_2$ (of $L$) and a cuspidal character $\delta \in \Irr(L_2)$.
Then $\eta$ is an irreducible constituent of $R^L_{L_2}(\delta)$. Then by the adjointness of the
Harish-Chandra induction and restriction and their transitivity \cite[Proposition 4.7]{DM}, we have that
$$\begin{aligned}0  & < c_\eta := [\SR^G_L(\chi),\eta]_{L} & = [\chi,R^G_L(\eta)]_G & \\
     & \leq [\chi,R^G_L(R^L_{L_2}(\delta))]_G
    & =     [\chi,R^G_{L_2}(\delta)]_G & \leq  [R^G_{L_1}(\psi),R^G_{L_2}(\delta)]_G.\end{aligned}$$
Since $\psi \in \Irr(L_1)$ and $\delta \in \Irr(L_2)$ are cuspidal,
it follows by \cite[Proposition 9.1.5]{C} that the pair $(L_1,\psi)$ is $G$-conjugate to the pair $(L_2,\delta)$ and
$R^G_{L_1}(\psi) = R^G_{L_2}(\delta)$. Hence, with no loss of generality we may replace  $(L_1,\psi)$ by
$(L_2,\delta)$. Furthermore, by \cite[Proposition 9.2.4]{C}, $[R^G_{L_1}(\psi),R^G_{L_1}(\psi)]_G$
can be bounded by the order of the Weyl group $W(\GC)$ of $\GC$ and so in terms of the semisimple rank $r$ as well.
Thus we can bound $c_\eta$ in terms of $r$. The same is true for
$[R^L_{L_1}(\psi),R^L_{L_1}(\psi)]_L$, and so for the number of possibilities for $\eta$.  In particular, if
$[\GC,\GC]$ is simple, then $|W(\LC)| \leq |W(\GC)| \leq W(r)$ and so we can choose $A(r) = W(r)^2$.
\hal

From now on we assume that $p$ is a good prime for $\GC$
(and $\KK = \overline\KK$ is a field of characteristic $p$).
 Then a theory of {\it generalized Gelfand-Graev representations}
(GGGRs) was developed by Kawanaka \cite{K}: for each unipotent element $u \in G = \GC^F$ one can associate a
GGGR with character $\G_u$ (which depends only the conjugacy class of $u$ in $G$).

Suppose now that $\OC=u^\GC$ is an
$F$-stable unipotent conjugacy class in $\GC$.
By the Lang-Steinberg theorem, since $\GC$ is connected we may assume
that $u \in G$. Then $\OC$ is called a {\it unipotent support} for a given $\rho \in \Irr(G)$ if
\begin{enumerate}[\rm(i)]
\item $\sum_{g \in \OC^F}\rho(g) \neq 0$;
\item If $\OC'$ is any $F$-stable unipotent class of $\GC$ such that $\sum_{g \in \OC'^F}\rho(g) \neq 0$, then
         $\dim \OC' \leq \dim \OC$.
\end{enumerate}
As shown in \cite{Geck}, as $p$ is a good prime for $\GC$, each $\rho \in \Irr(\GC^F)$ has a unique
unipotent support $\OC_\rho$ (see also \cite{GM}).

Next, $\OC \cap G$ is a disjoint union $\cup^r_{i=1}u_i^G$ of, say, $r$
conjugacy classes in $G$. If $A(x) = \CB_\GC(x)/\CB_\GC(x)^\circ$ is the component group of the centralizer of $x \in \GC$,
then one defines
$$\TG_u := \sum^r_{i=1}[A(u_i):A(u_i)^F]\G_{u_i}.$$
Then $\OC$ is called a {\it wave front set} for a given $\rho \in \Irr(G)$ if
\begin{enumerate}[\rm(i)]
\item $[\TG_u,\rho]_G \neq 0$;
\item If $\OC' = v^\GC$ is a unipotent class of $\GC$ with $v \in \GC^F$ such that $[\TG_v,\chi]_G \neq 0$, then
         $\dim\OC' \leq \dim\OC$.
\end{enumerate}
Work of Lusztig \cite{L} and subsequently \cite[Theorem 14.10]{Ta} show that each
$\rho \in \Irr(G)$ has a unique wave front set $\OCS_\rho$.
Moreover, if $\ZB(\GC)$ is connected, then $\OCS_\rho$ is the unipotent class denoted by
$\xi(\rho)$ in \cite[(13.4.3)]{L},
and, if $G$ is defined over $\F_q$, then as a polynomial in $q$ with rational coefficients, the degree of $\rho$ is
\begin{equation}\label{degree1}
  \rho(1) = \frac{1}{n_\rho}q^{(\dim \OCS_\rho)/2} + \mbox{ lower powers of }q,
\end{equation}
for some positive integer $n_\rho$ dividing $|A(u)|$ if $u \in \OC_\rho$.
Furthermore, if $D_G$ denotes the {\it Alvis-Curtis duality} (cf.
\cite[Chapter 8]{DM}), and $\rho^* = \pm D_G(\rho) \in \Irr(G)$ for $\rho \in \Irr(G)$, then
\begin{equation}\label{dual}
  \OC_{\rho^*} = \OCS_\rho,
\end{equation}
(see e.g. \cite[\S1.5]{Ta}).

\vspace{4mm}
The next two lemmas are well known to the experts. In particular, they have similar conclusions and proofs to Theorems 4.1(ii) and 1.7 of \cite{chardeg}. However, for application to bounding the function $f(r)$ in Theorem \ref{main1} (see Proposition \ref{frbd}), we need the extra detail in the lemmas concerning polynomials being products of
cyclotomic polynomials, which is not made explicit in \cite{chardeg}. We omit their proofs.

\begin{lem}\label{bounded2}
There is a constant $N = N(r)$ depending only on $r$ and a collection of $N$ monic polynomials, each being a product
of cyclotomic polynomials, such that the following statement holds.
If $\GC$ is a connected reductive group of semisimple rank $\leq r$ in good characteristic $p$,
$\GCF$ is defined over $\F_q$, and $s \in \GCF$ is semisimple, then
$$[\GC^F:(\CB_\GC(s)^\circ)^F]_{p'} = f(q),$$
where $f$ is one of the chosen polynomials.
\end{lem}

In what follows, with a slight abuse of language, we also view $t$ as a cyclotomic polynomial in variable $t$.

\begin{lem}\label{bounded3}
There are constants $B_1 = B_1(r)$ and $B_2 = B_2(r)$ depending only on $r$, and $B_2$ monic polynomials, each being
a product of cyclotomic polynomials in one variable $t$, such that the following statement
holds for any connected reductive algebraic group $\GC$ of semisimple rank $\leq r$ with connected center in
good characteristic.
When $\GCF$ is defined over $\F_q$ and $\chi \in \Irr(\GC^F)$, then
$$\chi(1) = (1/n_\chi)\Deg^*_\chi(q),$$
where $\Deg_\chi^*$ is one of the chosen monic polynomials,
$n_\chi \in \N$, $1 \leq n_\chi \leq B_1$. In fact, if $[\GC,\GC]$ is simple, then one can take $B_1$
to be the largest order of the component group $\CB_\HC(u)/\CB_\HC(u)^\circ$,
where $\HC$ is any simple algebraic group of rank $r$ and $u \in \HC$ any unipotent element.
\end{lem}

Recall that the set of unipotent classes in $\GC$ admit the partial order $\leq$, where
$u^\GC \leq v^\GC$ if and only if $u^\GC \subseteq \overline{v^\GC}$.

\begin{prop}\label{closure1}
Let $p$ be a good prime for $\GC$,  $G = \GC^F$, and let $u \in G$ be a unipotent element.
Then the following statements hold.

\begin{enumerate}[\rm(i)]
\item $D_G(\G_u)$ is unipotently supported, i.e. is zero on all non-unipotent elements of $G$.

\item Suppose that $D_G(\G_u)(v) \neq 0$ for some unipotent element $v \in G$. If
$\ZB(\GC)$ is disconnected, assume in addition that $q$ is large enough compared to
the semisimple rank of $\GC$. Then $u^\GC \leq v^\GC$.

\end{enumerate}
\end{prop}

\pf
(i) is well known, and (ii) is \cite[Scholium 2.3]{DLM2}. (Even though \cite{DLM2} assumes $p$ is large enough, in fact the
proof of \cite[Scholium 2.3]{DLM2} needs only that $p$ is a good prime.
As pointed out to the authors by J. Michel and J. Taylor, the proof in \cite{DLM2}
relies on the validity of the results in \cite{L3}, which were shown to hold under the indicated hypotheses
by Shoji \cite{Sh1}, cf. \cite[Theorem 4.2]{Sh2}.)
\hal


\begin{prop}\label{closure2}
Let $\GC/\ZB(\GC)$ be simple, $p$ be a good prime for $\GC$, $G = \GC^F$, and let $\ZB(\GC)$ be connected.
Suppose that $\chi \in \Irr(G)$ is such that $\SR^G_L(\chi) \neq 0$ for $L = \LC^F$,
where $\LC$ is a split Levi subgroup of $\GC$, and let
$\eta \in \Irr(L)$ be an irreducible constituent of $\SR^G_L(\chi)$. Let $\OCS_\chi = v^\GC$ and $\OCS_\eta = u^\LC$.
Then $\dim u^\GC \leq \dim v^\GC$.
\end{prop}

\pf
(i) To distinguish between GGGRs for $G$ and $L$, we will add the relevant superscript to their notation, e.g.
$\G^L_u$ is the GGGR of $L$ labeled by $u$.  First we show that if $R^G_L(D_L\G^L_u)(w) \neq 0$ for $w \in G$, then
$w$ is unipotent and $u^\GC \leq w^\GC$. Indeed, by Proposition \ref{closure1}(i), the generalized character
$D_L\G^L_u$ is unipotently supported, whence $R^G_L(D_L\G^L_u)$ is also unipotently supported. In particular,
$w$ is unipotent. Recall that $\LC$ is a Levi subgroup of an $F$-stable parabolic subgroup
$\PC$ with unipotent radical $\UC$. The condition on $w$ now implies that some $G$-conjugate of $w$ is $w' = xy$ where
$x \in \UC^F$, $y \in L$, and $D_L\G^L_u(y) \neq 0$. By
Proposition \ref{closure1} applied to $D_L\G^L_u$, $y$ is unipotent and $u^\LC \leq y^\LC$. It then follows
by \cite[Lemma 5.2]{GHM} (which is true for any connected reductive group $\GC$) that
$$u^\GC \leq y^\GC \leq (xy)^\GC = w^\GC,$$
as stated.

\smallskip
(ii) By the assumption, we may assume that $u \in L$ and $\eta$ is an irreducible constituent of the GGGR $\G^L_u$.
It follows that
$$0 < [\SR^G_L(\chi),\eta]_L \leq [\SR^G_L(\chi),\G^L_u]_L = [\chi,R^G_L(\G^L_u)]_G = $$
$$=   [D_G(\chi),D_G(R^G_L(\G^L_u))]_G = [D_G(\chi),R^G_L(D_L\G^L_u)]_G.$$
Here we use the self-adjointness of $D_G$ and the intertwining property of $D_G$ with
$R^G_L$ (see \cite[Proposition 8.10, Theorem 8.11]{DM}). In particular, there must exist some $w \in G$ such that
$$D_G(\chi)(w) \neq 0,~~R^G_L(D_L\G^L_u)(w) \neq 0.$$
Let $\chi^* = \pm D_G(\chi) \in \Irr(G)$ so that $\OC_{\chi^*} = \OCS_\chi = v^\GC$ with $v \in G$. By (i), the
condition $R^G_L(D_L\G^L_u)(w) \neq 0$ implies that $w$ is unipotent and
$$u^\GC \leq w^\GC.$$
Now we can apply \cite[Theorem 8.1]{AA} (which uses only
the assumption that $\ZB(\GC)$ is connected and $\GC/\ZB(\GC)$ is simple; cf. also
\cite[Corollary 13.6]{Ta}) to obtain from $\chi^*(w) \neq 0$ that
$$\dim w^\GC \leq \dim v^\GC.$$
It follows that
$$\dim u^\GC \leq \dim \overline{w^\GC} = \dim w^\GC  \leq \dim v^\GC,$$
as desired.
\hal

\vspace{4mm}
\no {\large {\bf Proof of Theorem \ref{main1}.} }

(i) Denoting $\rho = \SR^G_L(\chi)$, we have by Proposition \ref{fixed} that $|\chi(g)| = |\rho(g)| \leq \rho(1)$. Hence,
it suffices to bound $\rho(1)$ in terms of $\chi(1)$. Fix the semisimple rank $r$ of $\GC$.
First we handle the case where $\ZB(\GC)$ is connected. Note that $\HC := \GC/\ZB(\GC)$ is simple (of rank $r$)
as $[\GC,\GC]$ is  simple.
Consider any irreducible constituent $\eta$ of $\rho$ and let $\OCS_\eta = u^\LC$ for some $u \in L$ and
$\OCS_\chi = v^\GC$ for some $v \in G$. By Proposition \ref{closure2} we have $\dim u^\GC \leq \dim v^\GC$.
On the other hand, $\dim u^\LC \leq \a (\dim u^\GC)$ by the choice of $\a$, and so
\begin{equation}\label{degree2}
  \dim u^\LC \leq \a(\dim v^\GC).
\end{equation}
Now \eqref{degree1} and Lemma \ref{bounded3} imply that
$$\eta(1) \leq (q+1)^{(\dim u^\LC)/2},~~B_1\chi(1) \geq (q-1)^{(\dim v^\GC)/2}.$$
Let $D = D(r)$ denote the largest dimension of unipotent classes in simple algebraic groups of rank $r$.
Using \eqref{degree2} and noting that $\dim v^\GC = \dim v^\HC \leq D(r)$, we then get
$$\eta(1) \leq \left( \frac{q+1}{q-1} \right)^{\a D/2}\cdot B_1^\a\chi(1)^\a.$$
Setting $C := 3^{D/2}$ and applying Proposition \ref{hc}, we now obtain
$$\rho(1) \leq A(\max_\eta\eta(1)) \leq AB_1C\chi(1)^\a,
$$
and we are done in this case.

\smallskip
(ii) Next we handle the general case, where $\ZB(\GC)$ may be disconnected.
Consider a regular embedding of $\GC$ into $\tilde\GC$ with connected center and with compatible
Frobenius map $F:\tilde\GC \to \tilde\GC$, and set
$\tilde G := \tilde\GC^F$, $\ZC := \ZB(\tilde\GC)$. As $\tilde\GC = \ZC[\GC,\GC]$, $\tilde\GC$ and $\GC$ have the
same semisimple rank. Also, if $\LC$ is a Levi subgroup of an $F$-stable parabolic subgroup $\PC$ of $\GC$,
then we can embed $\PC$ in the $F$-stable parabolic subgroup $\tilde\PC = \UC\tilde\LC = \NB_{\tilde\GC}(\UC)$,
with the same unipotent radical $\UC$ as of $\PC$ and with $\tilde\LC = \ZC\LC$. Now, set $\tilde L:= \tilde\LC^F$ and note
that
\begin{equation}\label{levi}
  \tilde G = G\tilde L.
\end{equation}
Consider any $\chi \in \Irr(G)$ and some $\tilde\chi \in \Irr(\tilde G)$ lying above $\chi$, and denote
$$\rho := \SR^G_L(\chi),~~\tilde\rho := \SR^{\tilde G}_{\tilde L}(\tilde \chi).$$
Note that $\tilde\PC^F = U\tilde L$, and by \eqref{levi} we can
choose a set of representatives of $G$-cosets in $\tilde G$ that is contained in $\tilde L$. Hence, by Clifford's theorem
we can write
$$\tilde\chi|_G = \sum^t_{i=1}\chi^{x_i},$$
where $1=x_1, \ldots ,x_t \in \tilde L$. As $\tilde L$ normalizes $U$, we see that the Harish-Chandra restrictions $\rho_i$
of $\chi^{x_i}$ to the Levi subgroup $L$ all have the same dimension, equal to $[\chi|_U,1_U]_U$. Thus
$$\rho(1) = \tilde\rho(1)/t,~~\chi(1) = \tilde\chi(1)/t.$$
Now, any unipotent element $v \in \tilde L$ is contained in $\tilde L \cap \GC = L$, and
$v^{\tilde\GC} = v^\GC$ and similarly $v^{\tilde\LC} = v^{\ZC\LC} = v^\LC$. Thus the constants $\a$ for $L$ and for
$\tilde L$ as defined in Theorem \ref{main1} are the same. Applying Lemma \ref{fixed} to $\chi$ and the result of (i) to
$\tilde\chi$, we now have
$$|\chi(g)| = |\rho(g)| \leq \rho(1) = \frac{1}{t}\tilde\rho(1) \leq \frac{1}{t}f(r)\tilde\chi(1)^\a \leq f(r)\chi(1)^\a.$$
This completes the proof of Theorem \ref{main1}.
\hal

\vspace{2mm}
\no {\bf Remark } In the case of $GL_n(q)$, it is possible to give an alternate proof of Theorem \ref{main1} which does not use recent results on unipotent supports and
wave front sets; we do not give this here, but a sketch can be found in the last section of \cite{mwlprinceton}.

\vspace{4mm}
The next result provides a bound for the function $f$ in Theorem \ref{main1}.

\begin{prop}\label{frbd}
Under the assumptions of Theorem \ref{main1}, suppose that $q \geq q_0 \geq 2$.  Then $f(r)$ can be chosen to be
$$W(r)^2\cdot B(r) \cdot \left( \frac{q_0+1}{q_0-1} \right)^{(d(r)-r)/2},$$
where $W(r)$ is the largest order of the Weyl group of $\HC$, $B(r)$ is the largest order of $A(u)$ for unipotent
elements $u \in \HC$, and $d(r)$ is the largest dimension of $\HC$, when $\HC$ runs over simple algebraic groups
of rank $r$. In particular, if $r \geq 9$ and $q \geq r^2+1$, one can take
$$f(r) = 2^{2r+\sqrt{2r}+3}(r!)^2.$$
\end{prop}

\pf
By the proof of  Theorem \ref{main1} we may choose $f(r) =  AB_1C_1$, with $C_1 = ( \frac{q_0+1}{q_0-1})^{(d(r)-r)/2}$
(because $D = d(r)-r$). Next, $A \leq W(r)^2$ by Proposition \ref{hc} and $B_1 \leq B(r)$ by Lemma \ref{bounded3}.
Now assume that $r \geq 9$ and $q_0 \geq r^2+1$. Then $W(r) = 2^r \cdot r!$, $d(r) = 2r^2+r$ and so
$$\left( \frac{q_0+1}{q_0-1}\right)^{(d(r)-r)/2} \leq \left(1+\frac{2}{r^2}\right)^{r^2}.$$
It remains to bound $B(r)$. If $\HC = Spin_n$ (with $n = 2r$ or $2r+1$) and $u = \sum_i J_i^{r_i}$ is a unipotent element in $\HC$ with
$r_i$ Jordan blocks of size $i \geq 1$, then, according to \cite[\S3.3.5]{LS},
$|A(u)| \leq \max(2,2^k)$, where $k$ is the number of odd $i$ with $r_i > 0$. Note that
$$2r+1 \geq \sum^{k}_{j=1}(2j-1) = k^2,$$
and so $|A(u)| \leq 2^{\sqrt{2r+1}}$. Other simple groups of rank $r$ can be analyzed similarly using
\cite[Theorem 3.1]{LS} and yield smaller bound on $|A(u)|$. Hence we can take $B(r) = 2^{\sqrt{2r+1}}$ and complete
the proof by observing that
$$\left(1+\frac{2}{r^2}\right)^{r^2} \cdot 2^{\sqrt{2r+1}} < 2^{\sqrt{2r}+3}.$$
\hal

\vspace{4mm}
We conclude the section with some examples illustrating the sharpness of the $\a$-bound in Theorem \ref{main1}.

\begin{exa}\label{sharpness}
\begin{enumerate}[\rm(i)]

\item {\em  Let $G := GL_n(q)$ with $q > 2$, and let $g = \diag(\e,I_{n-1}) \in G$ for some $1 \neq \e \in \F_q^\times$.
Then $L := \CB_G(g) = GL_1(q) \times GL_{n-1}(q)$ is a Levi subgroup of $G$.

Let $\chi=\rho_n$ denote the unipotent character of $GL_n(q)$ labeled by the partition
$(n-1,1)$. Then $\rho_n(1) = (q^n-q)/(q-1)$. A computation inside the
Weyl group of $G$ (using the Comparison Theorem \cite[Theorem 5.9]{HL}) shows that
$$\SR^G_L(\rho_n) = 1_{GL_1(q)} \otimes (\rho_{n-1}+1_{GL_{n-1}(q)}).$$
Now Proposition \ref{fixed} implies that
$$|\chi(g)| = \rho_{n-1}(1)+1 = \frac{q^{n-1}-1}{q-1} \approx \chi(1)^{\frac{n-2}{n-1}}$$
if $q$ is large enough.
For this Levi subgroup $L$, the value of $\a$ in Theorem \ref{main1} is precisely $\frac{n-2}{n-1}$ (see Proposition \ref{gengl}), so the $\a$-bound is perfectly sharp in this example.

\item The Steinberg character $\St$ of a group $G = \GC^F$ as in Theorem \ref{main1} provides a good source of examples, since its values are easily calculated (see \cite[6.4.7]{C}): for a semisimple element $g \in G$,
\[
|\St (g) | = |\CB_G(g)|_p,
\]
(where $p$ is the underlying characteristic).

As a first example, let $G = GL_n(q)$ and let $g = \diag(\e,I_{n-1})$ as in the previous example. Then
\[
\St (g) = |GL_{n-1}(q)|_p = q^{\frac{1}{2}(n-1)(n-2)} = \St(1)^{\frac{n-2}{n}},
\]
while $\a = \frac{n-2}{n-1}$ for the Levi subgroup $\CB_G(g)$, as observed above.

As another example, let $G = GL_n(q)$ and suppose $n=mk$, where $2\le m\le q-1$ and $k>1$.
Let $\la_1,\ldots ,\la_m$ be distinct elements of $\F_q^\times$, and define
\[
g = {\rm diag}(\la_1 I_k,\cdots ,\la_m I_k) \in G.
\]
Let $L = \CB_G(g) = GL_k(q)^m$. By Corollary \ref{FL},
$\a(L) = \frac{1}{m}$.
On the other hand,
\[
\St (g) = q^{\frac{1}{2}mk(k-1)} = \St(1)^{\frac{k-1}{mk-1}},
\]
and the exponent $\frac{k-1}{mk-1}$ is close to $\a = \frac{1}{m}$ for $k$ large and $m$ fixed.

Similar examples showing the sharpness of Theorem \ref{main1} for the Steinberg character of other classical groups
can be constructed using \cite[Lemma 3.4]{lst2}.

\item Fix $m \geq 2$ and consider $G = GL_{2m}(q)$ with $q$ large enough (compared to $m$).
Again let $\la_1,\ldots ,\la_m$ be distinct elements of $\F_q^\times$, and define
\[
g = {\rm diag}(\la_1 I_2,\cdots ,\la_m I_2) \in G.
\]
Then $L = \CB_G(g) = GL_2(q)^m$, and $\al=\al(L) = 1/m$ as mentioned above. Consider the unipotent
characters $\chi^{(2m-j,j)}$ of $G$ labeled by the partition $(2m-j,j)$, $0 \leq j \leq m$. Then
$\sum^j_{i=0}\chi^{(2m-i,i)}$ is the permutation character of $G$ acting on the set of $j$-dimensional subspaces
of the natural module $V = \F_q^{2m}$. Note that $g$ fixes $mq^{m-1}(1+O(q^{-1}))$ $(m-1)$-dimensional subspaces of $V$, and
$(q+1)^m(1+O(q^{-1}))$ $m$-dimensional subspaces of $V$. It follows that for $\chi :=\chi^{(m,m)}$ we have
$$\chi(g) = q^m(1+O(q^{-1})),$$
whereas $\chi(1) = q^{m^2}(1+O(1/q))$. Thus $\chi(g) \approx \chi(1)^\al$.

\item More generally, fix $k,m \geq 2$ and consider $G = GL_{mk}(q)$ with $q$ large enough (compared to $\max(m,k)$).
Again let $\la_1,\ldots ,\la_m$ be distinct elements of $\F_q^\times$, and define
\[
g = {\rm diag}(\la_1 I_k,\cdots ,\la_m I_k) \in G.
\]
Then $L = \CB_G(g) = GL_k(q)^m$, and $\al=\al(L) = 1/m$ as mentioned above. Consider the unipotent
character $\chi:=\chi^\mu$ of $G$ labeled by the partition $\mu := (m^k) \vdash mk$. Observe that $\SR^G_L(\chi)$ contains
the Steinberg character $\St_L$ of $L$. (Indeed, by \cite[Proposition 5.3]{GLT}, the Alvis-Curtis duality functor $D_G$ sends $\chi$ to
$\pm \chi^\nu$, where $\nu = \mu'= (k^m) \vdash mk$, whereas $D_L(1_L) = \St_L$. Now, by \cite[Corollary 8.13]{DM} we have
$$\begin{array}{lll}
    [\SR^G_L(\chi),\St_L]_L & = [\SR^G_L(\pm D_G(\chi^\nu)),\St_L]_L & = \pm [D_L(\SR^G_L(\chi^\nu)),D_L(1_L)]_L\\
    & = \pm [\SR^G_L(\chi^\nu),1_L]_L & = \pm [\chi^\nu,R^G_L(1_L)]_G.\end{array}$$
But note that $L$ is a Levi subgroup of a rational parabolic subgroup of type $\nu$ of $GL_{mk}(\overline\F_q)$, whence
$\chi^\nu$ is an irreducible constituent of $R^G_L(1_L)$, and the claim follows.)  Since $\chi$ is a unipotent character and
the Harish-Chandra restriction preserves rational series, every irreducible constituent of $\SR^G_L(\chi)$ is a unipotent character of
$L$ and so contains $g \in \ZB(L)$ in its kernel. It now follows from Proposition \ref{fixed} that
$$\chi(g) = \SR^G_L(\chi)(g) = \SR^G_L(\chi)(1) \geq \St_L(1) = q^{mk(k-1)/2}.$$
On the other hand, the degree formula \cite[\S13.8]{C} implies that
$$\chi(1) = q^{m^2k(k-1)/2}(1+O(q^{-1})),$$
and we again obtain that $\chi(g) \gtrsim \chi(1)^\al$.

\item
As far as the exceptional groups of Lie type are concerned, it is again interesting to use the Steinberg character to test the sharpness of Theorem \ref{main1}. For example, let $G = E_8(q)$, and suppose $g \in G$ is a semisimple element with centralizer a Levi subgroup of type $E_7$. Then
\[
\St (g) = |E_7(q)|_p = q^{63} = \St (1)^\b,
\]
where $\b = \frac{21}{40}$, while the $\a$-value of this Levi is $\frac{17}{29}$, by Theorem \ref{alphaexcep}.
One can calculate such $\b$-values for all the Levi subgroups in Table \ref{extab} of Theorem \ref{alphaexcep}; it is never the case that $\b=\a$, but in some cases the values of $\b$ and $\a$ are quite close.

\item We offer one more example with $G = \GC^F = SL_n(q)$, with $q \geq n+2$, and
\[
g = \diag(\la_1,\la_2,\cdots ,\la_n) \in G,
\]
where $\la_1, \ldots,\la_n \in \F_q^\times$ are pairwise distinct.
Then $T=\CB_G(g)$ is a maximally split maximal torus.
Let $\mu \vdash n$ be such that the irreducible character $S^\mu$ of the Weyl group $W(\GC) \cong \SSS_n$ labeled by $\mu$
has the largest possible degree, and let $\chi:=\chi^\mu$ denote the unipotent character of $G$ labeled by $\mu$. As in (iv), every irreducible
constituent of $\SR^G_T(\chi)$ is trivial at $g$.
A computation in $W(\GC)$ and Proposition \ref{fixed} show that
$$\chi(g) = \SR^G_T(\chi)(g) = \SR^G_T(\chi)(1) = S^\mu(1),$$
whereas $\al(T) = 0$. Thus for the function $f$ in Theorem \ref{main1} we have
$$f(n-1) \geq S^\mu(1) \geq e^{-1.283\sqrt{n}}\sqrt{n!},$$
with the latter following from the main result of \cite{LoS}, \cite{VK}.
}
\end{enumerate}

\end{exa}

\section{General and special linear groups}\label{sec:main1b}

In this section we prove Theorems \ref{main2}, \ref{main1b} and
\ref{main1c}. Along the way we establish character bounds for unipotent elements of $GL_n(q)$ (see Theorem \ref{gl-uni}), and also for elements with extension-field centralizers for its semisimple parts
(Theorem \ref{gl-ext}).

\subsection{Proof of Theorem \ref{main1b}}\label{pf13}

We will keep the notation of \S\ref{sec:main1}.

\smallskip
(i) First we consider the case $\GC = GL_n$. In this case, the centralizer
of any element in $\GC$ is connected and one can check (e.g. using
\cite{Green}) that $n_\rho = 1$ in \eqref{degree1}.
Let $\varphi$ be an irreducible
$\ell$-Brauer character of $G = \GC^F = GL_n(q)$ and $g \in G$ as in Theorem \ref{main1b}. By Proposition
\ref{fixed}, $|\varphi(g)| = |\psi(g)|$ for $\psi := \SR^G_L(\varphi)$. According to \cite[Theorem B]{BK}, one can label
complex and $\ell$-Brauer characters of $G$ and find a complex character $\chi \in \Irr(G)$ with the same label as of
$\varphi$ such that both the generic degree of $\chi$ and the lower bound (given in \cite[Theorem B]{BK}) are
monic polynomials in $q$ of same degree say $N_\chi$. Using \eqref{degree1} and the equality $n_\chi = 1$, we have
$$N_\chi := (\dim \OCS_\chi)/2.$$
As $\chi(1)$ is a product of cyclotomic polynomials in $q$, we also have that
$$\chi(1) \leq (q+1)^{N_\chi}.$$
Furthermore, one can easily check that the lower bound in \cite[Theorem B]{BK} satisfies
$$\varphi(1) \geq (q-1)^{N_\chi}.$$
As $N_\chi \leq \dim \GC = n^2$, there is a constant $D=D(n) \leq 3^{n^2}$ such that
\begin{equation}\label{degree3}
  \frac{\chi(1)}{\varphi(1)} \leq D.
\end{equation}
Since $\chi$ and $\varphi$ have the same labeling, $\varphi$ is a constituent of the restriction $\chi^\circ$ of $\chi$
to $\ell'$-elements of $G$. Let $\PC = \UC\LC$ be an $F$-stable parabolic subgroup of $\GC$ with Levi subgroup $\LC$
and unipotent radical $\UC$.
As $U := \UC^F$ is an $\ell'$-group, we also have that
$$\psi(1) = [\varphi|_U,1_U]_U \leq [\chi|_U,1_U]_U = \rho(1),$$
where $\rho := \SR^G_L(\chi)$. The proof of Theorem \ref{main1} yields a function $f:\N \to \N$ such that
$$\rho(1) \leq f(n)\chi(1)^\a.$$
Choosing $h$ such that $h(n) \geq f(n)D$ and applying \eqref{degree3}, we obtain
$$|\varphi(g)| = |\psi(g)| \leq \psi(1) \leq \rho(1) \leq h(n)\varphi(1)^\a,$$
as desired.

\smallskip
(ii) Now we consider $S = SL_n(q)$ as the derived subgroup of $G = GL_n(q)$, an irreducible $\ell$-Brauer character
$\varphi_1$ of $S$, and its Harish-Chandra restriction $\psi_1$ to the Levi subgroup $L \cap S$, where
$L$ is a suitable split Levi subgroup of $G$. Note that we can
choose a set of representatives of $S$-cosets in $G$ that is contained in $L$. If $\varphi \in \IBR(G)$ lies above
$\varphi_1$, then by Clifford's theorem and the last observation we can write
$$\varphi|_S = \sum^t_{i=1}(\varphi_1)^{x_i},$$
where $1=x_1, \ldots ,x_t \in L$. As $L$ normalizes $U$, we see that the Harish-Chandra restrictions $\psi_i$
of $(\varphi_1)^{x_i}$ to the Levi subgroup $L \cap S$ all have the same dimension. Thus
$$\psi_1(1) = \psi(1)/t,~~\varphi_1(1) = \varphi(1)/t.$$
Applying Proposition \ref{fixed} and the results of (i), we now have
$$|\varphi_1(g)| = |\psi_1(g)| \leq \psi_1(1) = \frac{1}{t}\psi(1) \leq \frac{1}{t}h(n)\varphi(1)^\a \leq h(n)\varphi_1(1)^\a,$$
and so we are done.
\hal

\subsection{Proof of Theorem \ref{main2}}
(i) First we consider the case $G = GL_n(q)$. If $\LC$ is a torus, we can
choose a regular semisimple element $g \in L$ and take $\chi = 1_G$. Assume
now that $\LC$ is not a torus, and choose $u \in L_{\mathrm {unip}}$ such that
\begin{equation}\label{opt-1}
  \al(L) = (\dim u^\LC)/(\dim u^\GC).
\end{equation}
We may assume that
\begin{equation}\label{opt-1a}
L = GL_{n_1}(q) \times GL_{n_2}(q) \times \cdots \times GL_{n_r}(q),~~
  u = \diag(u_1,u_2, \ldots ,u_r)
\end{equation}
where $n_1 \geq n_2 \geq \ldots \geq n_r \geq 1$; furthermore, $u_i \in GL_{n_i}(q)$ is a
unipotent element, the sizes of whose Jordan blocks form a partition $\la_i \vdash n_i$.
Let $\mu_i \vdash n_i$ be the partition conjugate to $\la_i$ and let $\chi^{\mu_i}$
be the unipotent character of $GL_{n_i}(q)$ labeled by $\mu_i$. Now Green's formula for the degree
of $\chi^{\mu_i}$ (see the discussion before \cite[Theorem A]{BK}) implies that $\chi^{\mu_i}(1)$
is a monic polynomial in $q$ of degree $(1/2)\dim u_i^\LC$. Hence, if we choose
$C_n$ large enough, then using Lemma \ref{bounded2} we see that
\begin{equation}\label{opt-2}
  \gam(1) \geq (1/2)q^{(1/2)\dim u^\LC},
\end{equation}
for $\gam:= \chi^{\mu_1} \otimes \chi^{\mu_2} \otimes \ldots \otimes \chi^{\mu_r} \in \Irr(L)$
whenever $q \geq C_n$.

Next, let $\mu := \mu_1 + \mu_2 + \ldots + \mu_r$, where we have added zero parts to $\mu_i$
so that $\mu_1, \ldots ,\mu_r$ have the same number of parts, and then take the $i^{\mathrm {th}}$
part of $\mu$ to be the sum of all the $i^{\mathrm {th}}$ parts of $\mu_1, \ldots,\mu_r$. Again using
Green's formula, we then see that the unipotent character $\chi = \chi^\mu$ of $G$ labeled
by $\mu$ is a monic polynomial in $q$ of degree $(1/2)\dim^{u^\GC}$, whence
\begin{equation}\label{opt-3}
  \chi(1) \leq 2q^{(1/2)\dim u^\GC}
\end{equation}
if $q \geq C_n$.

For any $\nu \vdash m$, let $S^\nu$ denote the irreducible character of $\SSS_m$ labeled by $\nu$.
An application of the Littlewood-Richardson formula \cite[2.8.14]{JK} shows that the
restriction of $S^\mu$ to $\SSS_{n_1} \times \SSS_{n_2} \times \ldots \times \SSS_{n_r}$ contains
$S^{\mu_1} \otimes S^{\mu_2} \otimes \ldots \otimes S^{\mu_r}$. This computation in the Weyl groups of
$\GC$ and $\LC$ implies that $\gam$ is an irreducible constituent of $\SR^G_L(\chi)$. Now,
if $q \geq n+1$, then we can choose a semisimple element $g \in \ZB(L)$ such that
$\CB_G(g) = L$. As in Example \ref{sharpness}(iv), we have that every irreducible constituent of
$\SR^G_L(\chi)$ has $g$ in its kernel, and so by Proposition \ref{fixed},
$$\chi(g) = \SR^G_L(\chi)(g) = \SR^G_L(\chi)(1) \geq \gam(1).$$
Hence the statement follows from the choice \eqref{opt-1} of $u$ and the bounds \eqref{opt-2},
\eqref{opt-3}.

\smallskip
(ii) To handle the case of $SL_n(q)$, first we recall that
any unipotent character of $G$ remains irreducible over $SL_n(q)$. Furthermore,
as mentioned in the proof of Theorem \ref{main1}, for
any split Levi subgroup $\LC$ of $\GC = GL_n(\overline\F_q)$ we have $\GC = [\GC,\GC]\ZC$ and
$\LC = (\LC \cap [\GC,\GC])\ZC$, where $\ZC := \ZB(\GC)$. It follows that
$$\dim u^\GC = \dim u^{[\GC,\GC]},~~\dim u^\LC = \dim u^{\LC \cap [\GC,\GC]}$$
for any element $u \in \LC_{\mathrm {unip}} = (\LC \cap [\GC,\GC])_{\mathrm {unip}}$. Finally,
the condition $(n!)^n$ divides $q-1$ implies by the next Lemma \ref{det1} that, for any non-toral Levi subgroup $L$
given in \eqref{opt-1a} we can find an element $g \in SL_n(q)$ with $\CB_G(g) = L$. Now the
statement for $SL_n(q)$ follows from (i).
\hal

\begin{lem}\label{det1}
Let $1 \leq n_1 \leq n_2 \leq \ldots \leq n_r$ with $r,n_r \geq 2$, $n = \sum^r_{i=1}n_i$, and let $q$ be a prime power such that
$N:= n_r\cdot \prod^{r-1}_{i=1}(n_i+1)$ divides $q-1$. Then, for
$$L := GL_{n_1}(q) \times GL_{n_2}(q) \times \cdots \times GL_{n_r}(q) < G := GL_n(q),$$
there exists a semisimple element $s \in SL_n(q)$ such that $\CB_G(s) = L$.
\end{lem}

\pf
Choose $\zeta \in \F_q^\times$ of order $N$, and for any $d|N$ let $\zeta_d := \zeta^{N/d}$. Define
$$h_i = \zeta_{(n_1+1)(n_2+1) \ldots (n_i+1)}^{-1}I_{n_i} \in GL_{n_i}(q),~1 \leq i \leq r-1,~~
    h_r = \zeta I_{n_r} \in GL_{n_r}(q).$$
We prove by induction on $1 \leq i \leq r-1$ that $\prod^i_{j=1}\det(h_i) = \zeta_{(n_1+1)(n_2+1) \ldots (n_i+1)}$.
The induction base $i=1$ is obvious. For the induction step from $i-1$ to $i \geq 2$, we have
$$\prod^i_{j=1}\det(h_i) = \zeta_{(n_1+1)(n_2+1) \ldots (n_{i-1}+1)}\zeta_{(n_1+1)(n_2+1) \ldots (n_i+1)}^{-n_i} = \zeta_{(n_1+1)(n_2+1) \ldots (n_i+1)}.$$
Hence, for $s := \diag(h_1, h_2, \ldots,h_r) \in \ZB(L)$ we have
$$\det(s) =  \zeta_{(n_1+1)(n_2+1) \ldots (n_{r-1}+1)}\zeta^{-n_r} = 1.$$
The construction of $s$ and the condition $n_r \geq 2$ ensure that $\CB_G(s) = L$.
\hal

\subsection{Elements with extension-field centralizers}

\begin{thm}\label{gl-ext}
Let $G = GL_n(q)$ with $n \geq 2$ and $q \geq 8$, and let $g = su = us$ with $s \in G$ semisimple and $u \in G$ unipotent.
Suppose that $\CB_G(s) \cong GL_{n/k}(q^k)$ for some $1 < k \mid n$. Then
$$\frac{|\chi(g)|}{\chi(1)} \leq \frac{f(n)}{q^{n/3}}$$
for any $\chi \in \Irr(G)$ with $\chi(1) > 1$ and $f(n) = (11/7)^n-13/10$. In particular, if $q \geq 227$ then
$|\chi(g)| \leq \chi(1)^{1-1/2n}$ for all $\chi \in \Irr(G)$.
\end{thm}

\pf We proceed by induction on $n \geq 2$. Let $L:= SL_n(q)$ and let $W = \F_q^n$ denote the natural $G$-module.
Since $\chi(1) > 1$, all irreducible constituents of $\chi_L$
are non-trivial. In particular, if $n \geq 3$ then
$\chi(1) \geq (q^n-q)/(q-1) > q^{n-1}$ by \cite[Theorem 1.1]{TZ}.

\smallskip
(i) First we consider the case $k = n$. Since $|\CB_G(g)| \leq |\CB_G(s)| = q^n-1$, we have that
$|\chi(g)| \leq \sqrt{|\CB_G(g)|} < q^{n/2}$; in particular,
\begin{equation}\label{gl-ext-1}
  |\chi(g)|/\chi(1) \leq q^{1-n/2} \leq q^{-n/3}
\end{equation}
if $n \geq 6$, or if $n = 5$ and $\chi(1) \geq q^5$. The condition $\CB_G(s) \cong GL_1(q^n)$ also
implies that no eigenvalue of $g$ on $W$ can belong to $\F_q$.

Assume now that $n = 5$ and $\chi(1) < q^5$. By \cite[Theorem 3.1]{TZ}, every irreducible constituent of
$\chi_L$ is one of $q-1$ Weil characters, of degree $(q^n-1)/(q-1) - \delta$ with $\delta = 0$ or $1$. Since Weil
characters of $L$ extend to Weil characters of $G$, $\chi$ is a Weil character. Since no eigenvalue of $g$ on $W$
belongs to $\F_q$,  using the well-known character formula for
Weil characters of $G$, see e.g. \cite[(1.1)]{T}, we now see that $|\chi(g)| \leq q+1$ and so
$$|\chi(g)|/\chi(1) < (q^2-1)/(q^5-q) < q^{-n/3}.$$

Consider the case $n = 4$. If $\chi(1) \geq (q-1)(q^3-1)/2$, then
$$\frac{|\chi(g)|}{\chi(1)} \leq \frac{(q^4-1)^{1/2}}{(q^3-1)(q-1)/2} < q^{-n/3}$$
as $q \geq 8$. Assume now that  $\chi(1) < (q^3-1)(q-1)/2$. By \cite[Theorem 3.1]{TZ}, every irreducible constituent of
$\chi_L$ is one of $q-1$ Weil characters, all of which extend to Weil characters of $G$. Arguing as in
the previous case, we see that
$$|\chi(g)|/\chi(1) < (q^2-1)/(q^4-q) < q^{-n/3}.$$

If $n = 3$, then inspecting the character table of $G$ \cite{St} we get
$$\frac{|\chi(g)|}{\chi(1)} \leq \max\left(\frac{1}{q(q+1)},\frac{3}{(q^2-1)(q-1)}\right) < q^{-n/3}.$$
Similarly, for $n = 2$ we have $|\chi(g)|/\chi(1) \leq 2/(q-1) < (1.15)q^{-n/3}$ as $q \geq 8$.

Note that $f(n) > 1.16$ for all $n \geq 2$. Hence, to complete the induction base $2 \leq n \leq 5$,
it remains to consider the case $(n,k) = (4,2)$. Again inspecting the character table of $G$ \cite{Chevie}, we see that
$$|\chi(g)|/\chi(1) < 1/(q-1)^2 < q^{-n/3}.$$

\smallskip
(ii) From now on we may assume that $n \geq 6$ and $2 \leq k < n$. Consider the action of $u$ on the natural
module $W' = \F_{q^k}^{n/k}$ of $\CB_G(s)$. If this action induces an element with only one Jordan block, then
$|\CB_G(g)| = |\CB_{GL_{n/k}(q^k)}(u)| < q^n$ and again \eqref{gl-ext-1} holds. Thus we may assume that
the $\langle g\rangle$-module $W'$ is decomposable as a direct sum of two $\langle g\rangle$-submodules:
$W' = W'_1 \oplus W'_2$, with
$$\dim_{\F_{q^k}}W'_1 = a \geq n/2k,~\dim_{\F_{q^k}}W'_2 = b \geq 1.$$
Viewing $W'_i$ as vector spaces over $\F_q$, we get a $g$-invariant decomposition $W = W_1\oplus W_2$, with
$\dim W_1 = ak \geq n/2$ and $\dim W_2 = bk \geq 2$. Writing $g = \diag(g_1,g_2)$ with $g_i \in G_i:= GL(W_i)$ and
let $s_i$ denote the semisimple part of $g_i$, we have
$$\CB_{G_1}(s_1) = GL(W'_1) \cong GL_a(q^k),~~\CB_{G_2}(s_2) = GL(W'_2) \cong GL_b(q^k).$$
In particular, the induction hypothesis applies to the elements $g_i \in G_i$.

\smallskip
(iii) Let $V$ be a $\C G$-module affording the character $\chi$, and denote $L_i := [G_i,G_i]$ for
$i = 1,2$. We decompose the $G_1 \times G_2$-module $V$ as
$$V = V_1 \oplus V_2 \oplus V_3,$$
where $V_1 := \CB_V(L_1)$, every irreducible constituent of $V_2$ is trivial on $L_2$ but {\it not} on $L_1$, and
every irreducible constituent of $V_3$ is nontrivial on $L_1$ and on $L_2$. Let $\chi_j$ denote the
$G_1 \times G_2$-character afforded by $V_j$, for $1 \leq j \leq 3$.

If $\a \otimes \b$ is any irreducible constituent of $\chi_3$, then $\a(1), \b(1) > 1$ by the construction of $V_3$, whence
$$|\a(g_1)|/\a(1) \leq f(ak)q^{-ak/3},~~|\b(g_2)|/\b(1) \leq f(bk)q^{-bk/3}$$
by the induction hypothesis applied to $g_1 \in G_1$ and $g_2 \in G_2$. It follows that
\begin{equation}\label{gl-ext2}
  |\chi_3(g)|/\chi_3(1) \leq f(ak)f(bk)q^{-n/3}.
\end{equation}

Next, let $\a \otimes \b$ be any irreducible constituent of $\chi_2$. Then $\a(1) > 1$ and $\b(1) = 1$ by the construction of $V_2$, whence
$$|\a(g_1)|/\a(1) \leq f(a)q^{-ak/3},~~|\b(g_2)| = 1$$
by the induction hypothesis applied to $g_1 \in G_1$. It follows that
\begin{equation}\label{gl-ext3}
  |\chi_2(g)|/\chi_2(1) \leq f(a)q^{-ak/3}.
\end{equation}

\smallskip
(iv) We will now estimate $\chi_j(1)/\chi(1)$ for $j = 1,2$. Let $d(X)$ denote the smallest degree of a nontrivial
complex representation of a finite group $X$, and let
$$a_{m,q} := \frac{\sqrt{q-1}}{d(SL_m(q))} + \frac{1}{d(SL_{m-1}(q))},~~b_{m,q} := \sum^{\infty}_{i=m+1}a_{i,q}.$$
The proof of \cite[Proposition 4.2.3]{LST} shows that, if $U$ is any nontrivial irreducible $\C SL_n(q)$-module
for $n > m \geq 3$ and $SL_m(q)$ is embedded naturally in $SL_n(q)$, then
\begin{equation}\label{gl-ext4}
  \dim \CB_U(SL_m(q)) \leq b_{m,q}\dim U.
\end{equation}
By \cite[Theorem 1.1]{TZ}, for $m \geq 4$ we have
$a_{m,q} < (q+\sqrt{q-1})/q^{m-2}(q+1)$. It follows that
\begin{equation}\label{gl-ext5}
  b_{m,q} <\frac{q+\sqrt{q-1}}{q+1}\sum^{\infty}_{i=m-1}q^{-i} = \frac{q(q+\sqrt{q-1})}{q^{m-1}(q^2-1)} < \frac{1.36}{q^{m-1}}
\end{equation}
if $m \geq 3$ and $q \geq 8$. As $a_{3,q} = \sqrt{q-1}/(q^2+q) + e/(q-1)$ with $e := 3-\gcd(q,2)$ and $q \geq 8$, we then have
\begin{equation}\label{gl-ext6}
  b_{2,q} = a_{3,q} + b_{3,q} < 1.3q^{-2/3}.
\end{equation}
Now, since $ak \geq n/2 \geq 3$, we have $ak-1 \geq n/3$. Applying \eqref{gl-ext4} and \eqref{gl-ext5}, we get
\begin{equation}\label{gl-ext7}
  \chi_1(1)/\chi(1) \leq b_{ak}(q) < 1.36q^{-n/3}.
\end{equation}
Similarly,
\begin{equation}\label{gl-ext8}
  \chi_2(1)/\chi(1) \leq b_{bk}(q) < 1.36q^{-(bk-1)} \leq 0.17q^{-bk/3}.
\end{equation}
if $bk \geq 3$, and
\begin{equation}\label{gl-ext9}
  \chi_2(1)/\chi(1) \leq b_2(q) <  1.3q^{-2/3}.
\end{equation}
if $bk \geq 2$ (using \eqref{gl-ext6} instead of \eqref{gl-ext5}). Note that in the case $bk = 2$, we must have
$6 \leq n = ak+2$, and so $ak \geq 4$, $ak-1 \geq 1+n/3$, whence instead of \eqref{gl-ext7} we have
\begin{equation}\label{gl-ext10}
  \chi_1(1)/\chi(1) \leq b_{ak}(q) < 0.17q^{-n/3}.
\end{equation}

\smallskip
(v) Now, if $bk \geq 3$, then putting \eqref{gl-ext2}, \eqref{gl-ext3}, \eqref{gl-ext7}, \eqref{gl-ext8} together, we
obtain
$$|\chi(g)| \leq |\chi_1(1)| + |\chi_2(g)| + |\chi_3(g)| \leq \frac{\chi(1)}{q^{n/3}} \cdot (1.36+0.17f(ak)+f(ak)f(bk)).$$
If $bk = 2$, then \eqref{gl-ext2}, \eqref{gl-ext3}, \eqref{gl-ext9}, \eqref{gl-ext10} altogether imply
$$|\chi(g)| \leq |\chi_1(1)| + |\chi_2(g)| + |\chi_3(g)| \leq \frac{\chi(1)}{q^{n/3}} \cdot (0.17+1.3f(ak)+f(ak)f(bk)).$$
The choice $f(n) = (11/7)^n-1.3$ ensures that
$$f(n) = f(ak)f(bk) + 1.3f(ak)+1.3f(bk) + 0.39 > f(ak)f(bk) + 0.17f(ak)+ 1.36,$$
whence $|\chi(g)|/\chi(1) \leq f(n)q^{-n/3}$, completing the induction step of the proof. The last statement
then follows, since $f(n) < q^{n/12}$ when $q \geq 227$ and $\chi(1) < q^{n^2/2}$.
\hal


\subsection{Unipotent elements in general linear groups}

\begin{thm}\label{gl-uni}
There is a function $g:\N \to \N$ such that the following statement holds. For any $n \geq 2$, any prime
power $q$, $\ell = 0$ or any prime not dividing $q$, any irreducible $\ell$-Brauer character $\varphi$ of
$G:= GL_n(q)$, and any unipotent element $1 \neq u \in G$,
$$|\varphi(u)| \leq g(n) \cdot \varphi(1)^{\frac{n-2}{n-1}}.$$
\end{thm}

\pf
Note that the statement holds when $n=2$ (choosing $g(2) = 1$) as in this case we have $|\varphi(u)| \leq 1$. So in what follows
we may assume $n \geq 3$.

Recall the partial order $\leq$ on the set of unipotent classes of $\GC = GL_n(\overline\F_q)$:
$x^\GC \leq y^\GC$ precisely when $x^\GC \subseteq \overline{y^\GC}$, and we consider
$G = \GC^F$ for a suitable Frobenius endomorphism $F$.
Note that the unipotent classes in $\GC$ are parametrized by partitions of $n$. We will prove by induction using
the partial order $\leq$ that, if $u$ is parametrized by a partition $\lam \vdash n$ then
$$|\varphi(g)| \leq g_\lam(n) \cdot \varphi(1)^{\frac{n-2}{n-1}}$$
for some positive constant $g_\lambda(n)$ depending only on $\lam$. Then the statement follows by taking
$$g(n) := \max_{\lam \vdash n}g_\lam(n).$$

Observe that $u$ is a {\it Richardson} unipotent element, that is, we can find an $F$-stable parabolic subgroup
$\PC$ with unipotent radical $\UC$ such that $u^\GC \cap \UC$ is an open dense subset of $\UC$ that forms
a single $\PC$-orbit. Indeed, as shown in \cite[\S5.5]{Hu}, if $\mu = (\mu_1, \ldots ,\mu_l) \vdash n$ is the conjugate
partition associated to $\lam$, then one can just take $\PC$ to be the standard parabolic subgroup generated by
the upper-triangular Borel subgroup together with matrices in the block-diagonal form, with block sizes
$\mu_1 \times \mu_1, \ldots ,\mu_l \times \mu_l$.
Furthermore, $\CB_\GC(u)$ is connected (as $\GC = GL_n(\overline\F_q)$), of dimension
equal to $\dim\PC-\dim\UC$, and contained in $\PC$, see \cite[Corollary 5.2.2]{C}.
Since $u^\GC \cap \UC$ is a single $\PC$-orbit and $\CB_\GC(u) = \CB_\GC(u)^\circ = \CB_\PC(u)$ is connected,
by the Lang-Steinberg theorem, $u^\GC \cap \UC$ contains an $F$-stable element $u'$, i.e. $u' \in u^\GC \cap U$ for
$U := \UC^F$.
The connectedness of $\CB_\GC(u)$ implies by the Lang-Steinberg theorem that $u,u' \in u^\GC \cap U$ are
$G$-conjugate. Replacing $u$ by $u'$, we may assume that $u \in u^\GC \cap U$. Then, again applying the Lang-Steinberg
theorem, we see that any element $w \in u^\GC \cap U$ can be written as $huh^{-1}$ for some
$h \in P := \PC^F$. Conversely, the $P$-orbit of $u$ is contained in $u^\GC \cap U$. Thus $u^\GC \cap U$ is
a single $P$-orbit, and so
$$|u^\GC \cap U| = [P:C],$$
where $C := \CB_P(u) = \CB_\PC(u)^F = \CB_\GC(u)^F$. The structure of $\CB_\GC(u)$ is given in
\cite[Theorem 3.1]{LS}. As $\dim\CB_\GC(u) = \dim\PC-\dim\UC$, there is a constant $A(n)$ depending only on
$n$ such that
\begin{equation}\label{for-u}
  |u^\GC \cap U| \geq \frac{2}{3}|U|
\end{equation}
for all $q \geq A(n)$ and all $\lam \vdash n$. By taking $g(n)$ large enough, say
\begin{equation}\label{for-g}
  g(n) \geq \max_{q'= p^r < A(n)} \left\{ \frac{|\psi(w)|}{\psi(1)^{\frac{n-2}{n-1}}} \mid 1 \neq w \in GL_n(q'),~w \mbox{ unipotent},~\psi \in \IBR_\ell(GL_n(q'))
    \right\},
\end{equation}
we may assume that the condition $q \geq A(n)$ is indeed
satisfied.

Let $1 \neq v \in U \smallsetminus u^\GC$ be labeled by $\nu \vdash n$. Then
$$v \in \UC = \overline{u^\GC \cap \UC},$$
and so $v^\GC \leq u^\GC$. In particular, if $u^\GC$ is minimal with respect to $\leq$, then
no such $v$ exists. If $u^\GC$ is not minimal, then by the induction hypothesis applied to $v^\GC$ we have
\begin{equation}\label{for-v}
  |\varphi(v)| \leq g_\nu(n) \cdot \varphi(1)^{\frac{n-2}{n-1}}
\end{equation}
for some positive constant $g_\nu(n)$ depending only on $\nu$. We will let
$g'_\lam(n)$ be the largest among all $g_\nu(n)$ when $\nu$ runs over the partitions for all such $v$.

Let $\rho := \SR^G_L(\varphi)$, where $L$ is a Levi subgroup of $P$. Then
$$ \rho(1) = [\varphi_U,1_U]_U = \frac{1}{|U|}\left( \varphi(1) + \sum_{1 \neq v \in U \smallsetminus u^\GC}\varphi(v) +
      \sum_{w \in u^\GC \cap U}\varphi(w)\right),$$
and so
$$|u^\GC \cap U|\cdot|\varphi(u)|  \leq |U|\rho(1) + \sum_{1 \neq v \in U \smallsetminus u^\GC}|\varphi(v)| + \varphi(1).$$
It now follows from \eqref{for-u} and \eqref{for-v} that
$$|\varphi(u)| \leq \frac{3}{2}\rho(1) + \frac{1}{2}g'_\lam(n)\varphi(1)^{\frac{n-2}{n-1}} + \frac{3}{2|U|}\varphi(1).$$
The proof of Theorem \ref{main1b} and the bound $\al \le \frac{n-2}{n-1}$ in Proposition \ref{gengl} imply that
$$\rho(1) \leq h(n)\varphi(1)^{\frac{n-2}{n-1}}.$$
On the other hand,
$|U| \geq q^{n-1}$ and $\varphi(1) < q^{n^2/2}$, whence for $n \geq 4$ we have
$$ \frac{\varphi(1)}{|U|} < \varphi(1)^{\frac{n-2}{n-1}}.$$
The same conclusion holds for $n = 3$ since $\varphi(1) < q^4$ in this case.
Hence the statement follows for $u$ by taking
$$g_\lam(n) := \frac{3}{2}h(n) +  \frac{1}{2}g'_\lam(n) + \frac{3}{2}.$$
\hal

\subsection{Special linear groups}

\begin{prop}\label{gl-sl}
Let $n\geq 3$ and let $\F$ be an algebraically closed field of characteristic $\ell$, where either
$\ell = 0$ or $\ell \nmid q$. Let $V$ be an irreducible $\F GL_n(q)$-module which is reducible over
$SL_n(q)$. Then one of the following holds:
\begin{enumerate}[\rm(i)]
\item $\dim(V) > q^{(n^{2}+n)/4}(q-1)$.

\item $n = 3$ and $\dim(V) \geq (q-1)(q^2-1)$.

\item $n = 4$ and $\dim(V)  \geq (q-1)(q^{2}-1)(q^{3}-1)$.

\item $2|n$, and $\dim(V) = \prod^{n/2}_{j=1}(q^{2j-1}-1)$ or $\prod^{n/2}_{j=1}(q^{n/2+j}-1)/(q^j-1)$.
Furthermore, $V$ is as described in {\rm \cite[Proposition 5.10(ii), (iii)]{KT}}, and $V_{SL_n(q)}$ is a sum of two
irreducible constituents.
\end{enumerate}

\end{prop}

\pf
Repeat the same proof of \cite[Proposition 5.10]{KT}, but for all $n$.
\hal

\subsection{Proof of Theorem \ref{main1c}}

(i) In this proof, let $G := GL_n(q) = \GC^F$  in the notation of \S2, and let $S := [G,G]$.
Write $g = su = us$ with $s \in G$ semisimple and $u \in G$ unipotent.
By Theorem \ref{main1b} and Proposition \ref{gengl},
we may assume that there is no split Levi subgroup $\LC$ of $\GC$ such that $\CB_G(s) \leq \LC^F$,
equivalently, $\CB_G(s) \cong GL_{n/k}(q^k)$ for some $k|n$.
For a fixed $n$, by choosing $h(n)$ large enough (similarly to the choice \eqref{for-g} of $g(n)$),
we may assume that $q \geq 227$.
Hence, we are done by Theorems \ref{gl-ext} (when $k > 1$)
and \ref{gl-uni} (when $k = 1$) if $H = G$.

From now on we will assume that $H = S$, and let $\tchi \in \Irr(G)$ be lying above $\chi$. Applying the
result for $G$, we are done if $\chi=\tchi|_S$. Hence, we may assume that $\tchi|_S$ is reducible, and so either
\begin{equation}\label{sl-1}
  \chi(1) \geq \tchi(1)/[G:S] > q^{(n^2+n)/4}
\end{equation}
or case (iv) of Proposition \ref{gl-sl} holds for a $\C G$-module $V$ affording $\tchi$.
Since
$$\frac{1}{2}\prod^{n/2}_{j=1}(q^{2j-1}-1) > q^{n^2/4}(1-\sum^{n/2}_{j=1}q^{1-2j})/2 > q^{(n^2-1)/4}$$
when $2|n$ and $q \geq 227$, we now have that
\begin{equation}\label{sl-2}
  \chi(1) > q^{(n^2-1)/4}.
\end{equation}
Assume in addition that $g^S = g^G$.
By Clifford's theorem we may write $\tchi|_S = \sum^{t}_{i=1}\chi^{x_i}$ for some elements $x_i \in G$. Since
$g^S=g^G$, $g^{x_i}$ is $S$-conjugate to $g$ and so $\chi^{x_i}(g) = \chi(g)$. It follows that
$$|\chi(g)|/\chi(1) = |\sum^t_{i=1}\chi^{x_i}(g)|/t\chi(1) = |\tchi(g)|/\tchi(1),$$
and so we are done again. So we may assume that $g^S \neq g^G$.

\smallskip
(ii) Here we consider the case $k > 1$, and recall that $u$ is a unipotent element in $\CB_G(s) = GL_{n/k}(q^k)$.
Note that $\det_{\F_{q^k}}$ maps $\CB_G(s)$ onto $\F^\times_{q^k}$, and the norm map
$\F^\times_{q^k} \to \F_q^\times$ is surjective. It follows that $s^G = s^S$. Hence, our assumption $g^G \neq g^S$
implies that $u \neq 1$. It is well known that the centralizer of any non-central element in $GL_m(q)$ has order at most
$q^{m^2-2m+2}$. It follows that
$$|\CB_G(g)| = |\CB_{\CB_G(s)}(u)| \leq q^{n^2/k-2n+2k} \leq q^{n^2/2-2n+4},$$
whence $|\chi(g)| \leq |\CB_G(g)|^{1/2} \leq q^{n^2/4-n+2}$. Together with \eqref{sl-2}, this implies
$$|\chi(g)| < \chi(1)^{1-1/2n}.$$

\smallskip
(iii) Now we consider the case $k = 1$, i.e. $s \in \ZB(G)$, and prove the stronger bound that
\begin{equation}\label{bound-unip}
  |\chi(g)| \leq h(n)\chi(1)^{\frac{n-2}{n-1}}.
\end{equation}
 Without loss of generality we may assume that $g=u$.
Let $r_i$ denote the number of Jordan blocks of size $i$ in the Jordan canonical form of $u$ for each $i \geq 1$;
in particular, $\sum_i ir_i = n$. It is easy to see that $g^G = g^S$ if $\gcd(i \mid r_i\geq 1) = 1$.
So the assumption $g^G \neq g^S$ implies
\begin{equation}\label{sl-3}
  \gcd(i \mid r_i \geq 1) > 1,
\end{equation}
in particular, $r_1 = 0$.
We claim (for $n \geq 5$) that either
\begin{equation}\label{sl-4}
  |\CB_G(g)| \leq q^{(n^2-3n+6)/2}
\end{equation}
or $g$ has type $J_2^{n/2}$ (i.e. $r_2 = n/2$). Indeed, by \cite[Theorem 3.1]{LS} we have that
$|\CB_G(g)| < q^N$, where
\begin{equation}\label{for-N}
  N := \sum_i ir_i^2 + 2\sum_{i < j}ir_ir_j.
\end{equation}
Now, if $r_2 = 0$, then  $3N \leq (\sum_i ir_i)^2 = n^2$ and so \eqref{sl-4} holds for $n \geq 6$.
If $r_2 = 0$ and $n = 5$, then \eqref{sl-3} implies that $r_5 = 1$, again yielding \eqref{sl-4}.
Suppose now that $n/2 > r:= r_2  > 0$, whence $r_3 = 0$ by \eqref{sl-3} and
$n-2r = t:= \sum_{j \geq 4}jr_j \geq 4$. Then
$$N = 2r^2 + 4r\sum_{j \geq 4} r_j + \sum_{j \geq 4}jr_j^2 + 2 \sum_{4 \leq j < j'}jr_jr_{j'} \leq 2r^2 + rt + t^2/4
   \leq (n^2-4n+8)/2.$$
In the case of \eqref{sl-4}, $|\chi(g)| \leq q^{(n^2-3n+6)/4}$ and so \eqref{bound-unip} holds
because of \eqref{sl-2}.

It remains to consider the case $g = J_2^{n/2}$. Let $W = \F_q^n = \langle e_1, \ldots ,e_n \rangle_{\F_q}$
denote the natural module for $G$, with $g(e_1) = e_1$.  Here, $|\CB_G(g)| < q^{n^2/2}$ by \eqref{for-N}, whence
\begin{equation}\label{sl-5}
  |\chi(g)| \leq q^{n^2/4}.
\end{equation}
Suppose first that
\begin{equation}\label{sl-6}
  \chi(1) > q^{(n-1)(n-2)/2}.
\end{equation}
If $n \geq 8$, then \eqref{sl-5} and \eqref{sl-6} immediately imply \eqref{bound-unip}.
In the remaining case we have $n=6$. An application of Clifford's theorem to
the normal subgroup $SL_6(q)\ZB(GL_6(q))$ of $GL_6(q)$ yields
$2 \leq \tchi(1)/\chi(1) \leq 6$. In particular, in the case of Proposition \ref{gl-sl}(iv) we have
$$\chi(1) \leq \frac{(q^4-1)(q^5-1)(q^6-1)}{2(q-1)(q^2-1)(q^3-1)} <  q^{10},$$
contrary to \eqref{sl-6}. Thus Proposition \ref{gl-sl}(i) must hold, whence
$$\chi(1) \geq \frac{1}{6}q^{21/2}(q-1) \geq q^{45/4}$$
(for $q \geq 1301$, which can be guaranteed by taking $h(6)$ large enough). The latter,
together with \eqref{sl-5}, implies \eqref{bound-unip}.

It remains to consider the case where \eqref{sl-6} does not hold.
Let $\chi$ be afforded by a $\C S$-module $V$ and let $P := Stab_S(\langle e_1 \rangle_{\F_q}) = UL$. We
decompose the $P$-module $V$ as $\CB_V(U) \oplus [U,V]$ and let $\gam$, respectively $\delta$, denote the $P$-character of
$\CB_V(U)$, respectively of $[U,V]$. In particular, $\gam = \SR^S_L(\chi)$, and so, arguing as in part (ii) of the proof of
Theorem \ref{main1b} we get
\begin{equation}\label{sl-7}
  |\gam(g)| \leq \gam(1) \leq f(n)\chi(1)^{\frac{n-2}{n-1}}
\end{equation}
(for some function $f:\N \to \N$). Next, we decompose
$$[U,V] = \sum_{1_U \neq \la \in \Irr(U)}V_\la,$$
as a direct sum of $U$-eigenspaces, which are transitively permuted by $L \cong GL_{n-1}(q)$. Note that
$g$ has prime order $p|q$, and it acts on $\Irr(U) \smallsetminus \{1_U \}$ with exactly $q^{n/2}-1$ fixed points. Certainly, the trace of
$g$ in its action on $\sum_{\la \in \OC}V_\la$ for any nontrivial $g$-orbit $\OC$ on $\Irr(U) \smallsetminus \{1_U \}$ is zero.
Since $\chi(1) \leq q^{(n-1)(n-2)/2}$, we have that
$$|\delta(g)| \leq (q^{n/2}-1)\dim(V_\la) = (q^{n/2}-1) \cdot \frac{\dim([U,V])}{q^{n-1}-1} < \frac{\chi(1)}{q^{n/2-1}} \leq \chi(1)^{\frac{n-2}{n-1}}.$$
Together with \eqref{sl-7}, this completes the proof.
\hal

The above proof yields the following analogue of Theorem \ref{gl-uni}:

\begin{cor}\label{sl-uni}
Let $S := SL_n(q) \leq G:= GL_n(q)$, and let $u \in G$ be any nontrivial unipotent element.
Assume that either $\ell = 0$, or $\ell \nmid q$ and
$u^G = u^S$. Then for any $\varphi \in \IBR_\ell(S)$,
$$|\varphi(u)| \leq g(n) \cdot \varphi(1)^{\frac{n-2}{n-1}}.$$
\end{cor}


\begin{rem}
{\em For any $\varepsilon > 0$, it seems possible to improve the term $q^{n/3}$ in Theorem \ref{gl-ext} to $q^{n/(2+\varepsilon)}$ at the price
of using much bigger $f(n)$, as well as a much bigger lower bound on $q$. As a consequence, one could perhaps improve
the exponent $1-1/2n$ in Theorem \ref{main1c} to $1-1/((1+\varepsilon) n)$. But we did not try to pursue it.}
\end{rem}



\section{Bounds for the constant $\a(\LC)$: Proof of Theorems~\ref{ratio}, \ref{alphaexcep}
and \ref{GL}}

For the proof of Theorem \ref{ratio}, it is convenient to handle the classical types $SL$, $Sp$ and $SO$ separately.
As in the theorem, let $\KK$ be an algebraically closed field of good characteristic. Note that by the defnition of $\a(\LC)$, this value does not depend on the isogeny type of $\GC$.

\subsection{Case $\GC = GL_n(\KK)$ or $SL_n(\KK)$}

To prove Theorem \ref{ratio} in this case we use the following lemma, which transfers attention from unipotent to semisimple elements in the analysis of $\a(\LC)$. Denote by $J_i$ a unipotent $i\times i$ Jordan block matrix, and by $\sum_i J_i^{n_i}$ the matrix in $SL_n(\KK)$ with $n_i$ diagonal blocks $J_i$ for each $i$, where $n = \sum in_i$.

\begin{lem}\label{uss}
Let $u = \sum_iJ_i^{n_i}$ be a unipotent element of $GL_n(\KK)$ where $n = \sum in_i$. Let $\la_j~(j\in \N)$ be distinct scalars in $\KK^\times$,
and for each $i$ let $D_i = {\rm diag}(\la_1,\la_2,\ldots ,\la_i)$. Define
\[
s := \sum_i D_i^{n_i} \in GL_n(\KK).
\]
Then $\dim \CB_{GL_n(\KK)}(u) = \dim \CB_{GL_n(\KK)}(s)$.
\end{lem}

\pf Observe that
\[
\begin{array}{ll}
\dim \CB_{GL_n(\KK)}(s) & = \left(\sum_{i\ge 1} n_i\right)^2 +  \left(\sum_{i\ge 2} n_i\right)^2 + \cdots \\
                                     & = \sum_i in_i^2 + 2\sum_{i<j}in_in_j
\end{array}
\]
which is equal to $\dim \CB_{GL_n(\KK)}(u)$ by \cite[3.1]{LS}. \hal

\vspace{4mm}
For a subgroup $X$ of $GL_n(\KK)$, define $X_{\mathrm {ss}} $ to be the set of semisimple elements of $X$.

\begin{cor}\label{sscomp}
If $n \geq 2$ and $\LC$ is a Levi subgroup of $\GC = GL_n(\KK)$ or $SL_n(\KK)$, then
\[
\a(\LC) \leq \hbox{{\rm {max}}}_{s \in \LC_{\mathrm {ss}} \smallsetminus \ZB(\GC)} \frac{\dim s^{\LC}}{\dim s^{\GC}}.
\]
\end{cor}

\pf We have $\LC = \GC \cap \prod_{j=1}^rGL_{a_j}(\KK)$, where $\sum_{j=1}^r a_j=n$. Let $u \in \LC_{\mathrm {unip}}$, so that
\[
u = \sum_{j=1}^r \sum_iJ_i^{n_{ij}},
\]
where $\sum_i in_{ij} = a_j$. The condition $u \neq 1$ means that there are some $i \geq 2$ and $j \geq 1$ such that
$n_{ij} > 0$. If we define $s = \sum_{j=1}^r \sum_iD_i^{n_{ij}}$, where $D_i$ is as in the statement of Lemma \ref{uss} (and the scalars $\la_j$ are chosen so that $s$ has determinant 1 in the case where $\GC = SL_n(\KK)$), then $s \notin \ZB(\GC)$.
Now Lemma \ref{uss} shows that $\dim \CB_{\LC}(u) = \dim \CB_{\LC}(s)$ and $\dim \CB_{\GC}(u) = \dim \CB_{\GC}(s)$.
\hal

\vspace{4mm}
\no {\bf Proof of Theorem \ref{ratio} for $GL_n(\KK)$, $SL_n(\KK)$}

\vspace{2mm}
We prove the theorem for $\GC = GL_n(\KK)$ and point out the small adjustment needed for $SL_n(\KK)$ at the end of the proof. Let $\LC$ be a Levi subgroup of $\GC$. Adopting an obvious
notational convention we take
\[
\LC = GL_a(\KK) \times GL_b(\KK) \times \cdots \times GL_z(\KK).
\]
Write $V_n = V_a \oplus V_b \oplus \cdots \oplus V_z$ for the corresponding
direct sum decomposition of $V_n = \KK^n$.


In view of Corollary \ref{sscomp}, it suffices to prove that
\begin{equation}\label{ineqa}
\hbox{max}_{s \in \LC_{\mathrm {ss}} \smallsetminus \ZB(\GC)} \frac{\dim s^{\LC}}{\dim s^{\GC}} \le \frac{1}{2}\left(1+ \frac{\dim \LC}{\dim \GC}\right).
\end{equation}

Let $s$ be a semisimple element of $\LC$, and let $\la_1,\ldots ,\la_k$ be
the distinct eigenvalues of $s$ on $V_n$. Write
\[
s|_{V_a} = \hbox{diag}\,(\la_1^{(a_1)},\ldots ,\la_k^{(a_k)}),
\ldots ,s|_{V_z} = \hbox{diag}\,(\la_1^{(z_1)},\ldots ,\la_k^{(z_k)}),
\]
where $\sum_{i=1}^k a_i = a$, and so on (superscripts denote
multiplicities). Then
\[
\CB_\GC(s) = GL_{a_1+b_1+\ldots }(\KK) \times \cdots \times GL_{a_k+b_k+\ldots }(\KK),
\]
\[
\CB_\LC(s) = \prod_{i=1}^k GL_{a_i}(\KK) \times \prod_{i=1}^k GL_{b_i}(\KK) \times \cdots .
\]
To prove (\ref{ineqa}) we need to show
\begin{equation}\label{try}
{1\over 2}{{\dim \GC  - \dim \LC} \over {\dim \GC}} \le {{\dim s^\GC - \dim s^\LC}
\over {\dim s^\GC}}.
\end{equation}
Now
\[
{1\over 2}{{\dim \GC  - \dim \LC} \over {\dim \GC}} = {{ab+ac+bc+\ldots } \over
{(a+b+\ldots )^2}},
\]
while
\[
{{\dim s^\GC - \dim s^\LC} \over {\dim s^\GC}} =
{{\sum_{i\ne j} (a_ib_j+a_ic_j+b_ic_j+\ldots )} \over
{\sum_{i<j} (a_i+b_i+\cdots )(a_j+b_j+\cdots )}}.
\]
Hence (\ref{try}) is equivalent to the inequality
\begin{equation}\label{ineqb}
\begin{array}{l}
\left(\sum a_i\sum b_i + \sum a_i \sum c_i + \sum b_i \sum c_i + \cdots\right)
\cdot \left(\sum_{i<j} (a_i+b_i+\cdots )(a_j+b_j+\cdots )\right) \\
\le \left(\sum a_i + \sum b_i + \cdots \right)^2 \cdot \sum_{i\ne j} (a_ib_j + a_ic_j
+ b_ic_j+ \cdots ).
\end{array}
\end{equation}
Now observe that all the terms on the left hand side of this inequality appear
with at most the same multiplicity on the right hand side. Hence (\ref{ineqb})
holds, and the proof is complete for $\GC = GL_n(\KK)$.

For the case where $\GC = SL_n(\KK)$, we need to prove the inequality (\ref{ineqb}) with the first term on the right hand side replaced by $\left(\sum a_i + \sum b_i + \cdots \right)^2-1$. This remains true, since the terms on the right hand side but not the left hand side of (\ref{ineqb}) include $\sum_{i\ne j} (a_i^3b_j + a_i^3c_j + b_i^3c_j+ \cdots )$, which is at least $\sum_{i\ne j} (a_ib_j + a_ic_j + b_ic_j+ \cdots ) \geq 1$. \hal

\vspace{4mm}
We also deduce the following general bound, which was used in the proof of Theorem \ref{gl-uni} and also in Example \ref{sharpness}.

\begin{prop}\label{gengl}
If $\LC$ is a Levi subgroup of $\GC = GL_n(\KK)$, then $\a(\LC) \le \frac{n-2}{n-1}$, with equality if and only if $\LC = GL_{n-1}(\KK) \times GL_1(\KK)$.
\end{prop}

\pf Choose maximal $a\le \frac{n}{2}$ such that $\LC \le GL_a(\KK)\times GL_{n-a}(\KK)$. By Theorem \ref{ratio} for $GL_n(\KK)$, proved above, we have
\[
\a(\LC) \le \frac{1}{2}\left(1+\frac{\dim \LC}{\dim \GC}\right) \le \frac{1}{2}\left( 1+\frac{a^2+(n-a)^2}{n^2}\right).
\]
One checks that the right hand side is less than $\frac{n-2}{n-1}$ for $n\ge 2a$, except in the following cases:
\begin{itemize}
\item[(a)] $a=1$, in which case $\LC = GL_1(\KK)\times GL_{n-1}(\KK)$ (by the maximal choice of $a$);
\item[(b)] $a=2$, $n\le 5$.
\end{itemize}
In case (b) we compute the values of $\a(\LC)$ and find that $\a(\LC) \le \frac{1}{2} < \frac{n-2}{n-1}$ (note that $n \ge 2a = 4$ in this case).

Hence it remains to consider case (a). We claim that in this case, $\a(\LC) = \frac{n-2}{n-1}$. Let $u$ be a nontrivial unipotent element of $\LC$, and write $u = \sum J_i^{n_i}$, where $\sum in_i = n$. Then $u$ projects to the element
$J_1^{n_1-1}+\sum_{i\ge 2}J_i^{n_i}$ in the factor $GL_{n-1}(\KK)$ of $\LC$, so by \cite[3.1]{LS}, we have
\[
\begin{array}{l}
\dim \CB_{\GC}(u) = \sum in_i^2 + 2\sum_{i<j}in_in_j, \\
\dim \CB_{\LC}(u) = 1+(n_1-1)^2 +  \sum_{i\ge 2} in_i^2 + 2(n_1-1)\sum_{j\ge 2}n_j + 2\sum_{2\le i<j}in_in_j.
\end{array}
\]
Defining $s:= \dim [V,u] = n-\sum n_i$, it follows that
\[
(\dim \GC - \dim \LC)-(\dim \CB_\GC (u)-\dim \CB_\LC (u)) = 2s.
\]
Next observe that
\[
\begin{array}{ll}
\frac{\dim u^\LC}{\dim u^\GC} \le \frac{n-2}{n-1}& \Leftrightarrow (n-1)\left((\dim \GC - \dim \LC)-(\dim \CB_\GC (u)-\dim \CB_\LC (u))\right) \ge \dim u^\GC  \\
& \Leftrightarrow 2(n-1)s \ge \dim u^\GC.
\end{array}
\]
By \cite[3.4(i)]{lish99} and its proof, we have $\dim u^\GC \le s(2n-s)$, so the above inequality holds when $s \ge 2$.
Finally, when $s=1$ we have $u = J_2+J_1^{n-2}$, and we calculate that $\frac{\dim u^\LC}{\dim u^\GC}
=  \frac{n-2}{n-1}$. Hence $\a(\LC) = \frac{n-2}{n-1}$ in case (1), and the proof is complete. \hal

\subsection{Symplectic groups}

Now we prove Theorem \ref{ratio} for symplectic groups. We revert to  Lie-theoretic notation, so assume that
\[
\GC = C_n = Sp_{2n}(\KK)  = Sp(V),
\]
where $V = V_{2n}(K)$ is the natural module for $\GC$ and $n\ge 2$.

Let $\LC$ be a Levi subgroup of $\GC$, so that $\LC' =  C_{n-r}\times \prod A_{r_i} \le C_{n-r} \times A_{r-1} $, where $1\le r\le n$. The first lemma deals with the case where $r=n$.

\begin{lem}\label{firstcn}
If $\LC \le A_{n-1}T_1$, then $\a(\LC) \le \frac{1}{2}$.
\end{lem}

\pf Assume $\LC \le A_{n-1}T_1 = GL_n$, and let $u$ be a nontrivial unipotent element of $\LC$. Write $u = \sum J_i^{n_i} \in SL_n$, where $\sum in_i = n$. As an element of $G = Sp_{2n}$, $u$ has Jordan form $\sum J_i^{2n_i}$. Hence by \cite[3.1]{LS},
\[
\begin{array}{l}
\dim \CB_{GL_n}(u) = \sum in_i^2+2\sum_{i<j}in_in_j := c_u, \\
\dim \CB_{\GC}(u) = 2\sum in_i^2+4\sum_{i<j}in_in_j + \sum_{i\;odd}n_i.
\end{array}
\]
So $\dim \CB_{\GC}(u) = 2c_u + s_u$, where $s_u = \sum_{i\;odd}n_i$. It follows that
\[
\frac{\dim u^\LC}{\dim u^\GC} \le \frac{\dim u^{GL_n}}{\dim u^\GC}\le \frac{n^2-c_u}{2n^2+n-2c_u-s_u} \le \frac{1}{2},
\]
and the conclusion follows. \hal

\begin{lem}\label{secondcn}
If $\LC = C_{n-r}T_r$, then $\a(\LC) \le \frac{1}{2}\left( 1+ \frac{\dim \LC}{\dim \GC}\right)$.
\end{lem}

\pf Let $u$ be a nontrivial unipotent element of $\LC' = C_{n-r} = Sp_{2n-2r}$, and write $u = \sum J_i^{n_i}$ with $\sum in_i = 2n-2r$. In $G = Sp_{2n}$, $u$ has Jordan form $J_1^{n_1+2r}+\sum_{i\ge 2}J_i^{n_i}$. Using \cite[3.1]{LS}, we find that
\begin{equation}\label{diff}
\dim \CB_{\LC}(u) - \dim \CB_{\GC}(u) = 2r\sum n_i + 2r^2.
\end{equation}
As in \cite[p.509]{lish99}, define
\[
s : = \dim [V,u] = 2n-2r - \sum n_i.
\]
Then (\ref{diff}) implies that $\dim u^\GC - \dim u^\LC = 2rs$. It also follows from \cite[3.4]{lish99} and its proof that
\begin{equation}\label{99res}
\dim u^\GC \le \frac{1}{2}s(4n-s+1).
\end{equation}
Now observe that
\[
\begin{array}{ll}
\frac{\dim u^\LC}{\dim u^\GC} \le \frac{1}{2}\left( 1+ \frac{\dim \LC}{\dim \GC}\right)
& \Leftrightarrow \dim u^\GC \le \frac{2\dim \GC ( \dim u^\GC - u^\LC)}{\dim \GC - \dim \LC} \\
&  \Leftrightarrow \dim u^\GC \le \frac{2s(2n^2+n)}{2n-r}.
\end{array}
\]
Clearly (\ref{99res}) implies that the last inequality holds, and so we are done. \hal

\begin{lem}\label{thirdcn}
Suppose $\LC' \le C_{n-r}\times A_{r-1} \le C_{n-r}\times C_r < \GC$ with $r>1$, and let $u = u_1u_2 \in \LC$ be a unipotent element with $u_1 \in C_{n-r}$, $u_2 \in A_{r-1} < C_r$. Then
\[
\dim u^\GC \ge \dim u_1^\GC + \dim u_2^{C_r}.
\]
\end{lem}

\pf Let $u_1 = \sum J_i^{a_i} \in C_{n-r}$ and $u_2 = \sum J_i^{b_i} \in A_{r-1}$, where $\sum ia_i = 2n-2r$, $\sum ib_i = r$. Then
\[
\begin{array}{l}
u_1 = J_1^{a_1+2r} + \sum_{i\ge 2}J_i^{a_i} \in \GC = Sp_{2n}, \\
u_2 = \sum J_i^{2b_i} \in C_r, \hbox{ and }\\
u = \sum J_i^{a_i+2b_i} \in \GC.
\end{array}
\]
Now using \cite[3.1]{LS}, we compute that
\[
\begin{array}{ll}
& \dim \CB_\GC(u_1)+\dim \CB_{C_r}(u_2)-\dim \CB_\GC(u) \\
= & 2r^2+r+2r\sum a_i - 2\sum ia_ib_i -2\sum_{i<j}i(a_ib_j+a_jb_i) \\
 =  &  2r^2+r + 2\left((\sum a_i)(\sum ib_i) - \sum ia_ib_i -\sum_{i<j}i(a_ib_j+a_jb_i)\right) \\
 \ge &  2r^2+r \\
  = &  \dim C_r,
\end{array}
\]
and the result follows. \hal

\vspace{4mm}
\no {\bf Proof of Theorem \ref{ratio} for $\GC = C_n$}

\vspace{2mm}
Let $\LC$ be a Levi subgroup of $\GC$, so $\LC' = C_{n-r}\times \prod A_{r_i} \le C_{n-r} \times A_{r-1} \le C_{n-r} \times C_r$, where $1\le r\le n$. Let $u = u_1u_2$ be a nontrivial unipotent element of $\LC$, where $u_1 \in C_{n-r}$, $u_2 \in A_{r-1}$. Using Lemma \ref{thirdcn}, we have
\begin{equation}\label{finalcn}
\frac{\dim u^\LC}{\dim u^\GC} \le \frac{\dim u_1^{C_{n-r}} + \dim u_2^{A_{r-1}}}{\dim u_1^\GC + \dim u_2^{C_r}}.
\end{equation}
Also Lemmas \ref{firstcn} and \ref{secondcn} imply that
\[
\frac{\dim u^{A_{r-1}}}{\dim u^{C_r}} \le \frac{1}{2},\;\; \hbox{ and } \frac{\dim u^{C_{n-r}}}{\dim u^\GC} \le \frac{1}{2}\left(1+
\frac{\dim C_{n-r}T_r}{\dim \GC}\right).
\]
Hence (\ref{finalcn}) implies that
\[
\frac{\dim u^\LC}{\dim u^\GC} \le  \frac{1}{2}\left(1+ \frac{\dim C_{n-r}T_r}{\dim \GC}\right)
\le  \frac{1}{2}\left(1+ \frac{\dim \LC}{\dim \GC}\right).
\]
This completes the proof of Theorem \ref{ratio} for $\GC = C_n$.

\subsection{Orthogonal groups}

We complete the proof of Theorem \ref{ratio} by handling the orthogonal groups. The proof for $\GC = B_n = SO_{2n+1}(\KK)$ is very similar to that for $\GC = C_n$: one shows that Lemmas \ref{firstcn}--\ref{thirdcn} also hold in the $B_n$ case (with $\LC = B_{n-r}T_r$ in Lemma \ref{secondcn} and $\LC'\le B_{n-r}\times A_{r-1} \le B_{n-r}\times D_r$ in Lemma \ref{thirdcn}), and the theorem follows. Things are a little different in the $D_n$ case, so assume now that
\[
\GC = D_n = SO_{2n}(\KK) = SO(V) \;\;(n\ge 4).
\]

Let $\LC$ be a Levi subgroup of $\GC$. Then $\LC' = D_{n-r} \times \prod A_{r_i} \le D_{n-r} \times A_{r-1}$, where $1 \le r \le n$ and $r \ne n-1$.

\begin{lem}\label{firstdn}
Suppose $\LC' \le A_{n-1}$ and $\LC'\ne A_{n-1}$. Then $\a(\LC) \le \frac{1}{2}$.
\end{lem}

\pf By assumption, $\LC \le GL_a \times GL_b$ where $a+b = n$ and $a,b\ge 1$. Let $u = u_1u_2 \in  \LC$, where
$u_1 = \sum J_i^{a_i} \in GL_a$ and $u_1 = \sum J_i^{b_i} \in GL_b$ (so $\sum ia_i = a$, $\sum ib_i = b$). Then $u = \sum J_i^{2a_i+2b_i} \in \GC$. By \cite[3.1]{LS},
\[
\dim \CB_\LC(u) = \sum ia_i^2 + \sum ib_i^2 + 2\sum ia_ia_j + 2\sum ib_ib_j =: c_u,
\]
and
\[
\begin{array}{ll}
\dim \CB_\GC(u) & = 2\sum i(a_i+b_i)^2 + 4\sum_{i<j} i(a_i+b_i)(a_j+b_j) - \sum_{i\;odd}(a_i+b_i) \\
                        & = 2c_u + 4\sum ia_ib_i + 4\sum_{i<j}i(a_ib_j+a_jb_i) - \sum_{i\;odd}(a_i+b_i).
\end{array}
\]
Then $\dim u^\LC = a^2+b^2-c_u$, while
\[
\dim u^\GC = 2(a+b)^2-(a+b)-(2c_u + 4\sum ia_ib_i + 4\sum_{i<j}i(a_ib_j+a_jb_i) - \sum_{i\;odd}(a_i+b_i)).
\]
To prove the lemma we need to show that $\dim u^\GC \ge 2(a^2+b^2-c_u)$. Using the equations $\sum ia_i = a$, $\sum ib_i = b$, this amounts to showing that
\begin{equation}\label{ineqz}
4\sum (i^2-i)a_ib_i + 4\sum_{i<j}i(j-1)(a_ib_j+a_jb_i)+ \sum_{i\;odd}(a_i+b_i) \ge \sum ia_i+\sum jb_j.
\end{equation}
Consider a term $ka_k+lb_l$ on the right hand side, with $a_k,b_l \ne 0$. If $k=l=1$ this occurs in the sum $\sum_{i\;odd}(a_i+b_i)$; if $k=l\ge 2$ it is less then or equal to the term $4 (k^2-k)a_kb_k$ on the left hand side; and if $k<l$ or $l<k$, it is at most $4k(l-1)a_kb_l$ or $4l(k-1)a_kb_l$, respectively. Hence the inequality (\ref{ineqz}) holds, completing the proof of the lemma. \hal

The proofs of the next two lemmas are very similar to those of Lemmas \ref{secondcn} and \ref{thirdcn}.

\begin{lem}\label{seconddn}
If $\LC = D_{n-r}T_r$, then $\a(\LC) \le \frac{1}{2}\left( 1+ \frac{\dim \LC}{\dim \GC}\right)$.
\end{lem}

\begin{lem}\label{thirddn}
Suppose $\LC' \le D_{n-r}\times A_{r-1} \le D_{n-r}\times D_r < \GC$ with $r>1$, and let $u = u_1u_2 \in \LC$ be a unipotent element with $u_1 \in D_{n-r}$, $u_2 \in A_{r-1} < D_r$. Then
\[
\dim u^\GC \ge \dim u_1^\GC + \dim u_2^{D_r}.
\]
\end{lem}

Let $\LC' = D_{n-r} \times \LC_1 \le D_{n-r} \times A_{r-1}$, where $\LC_1 = \prod A_{r_i} \le A_{r-1}$. If either $\LC_1 < A_{r-1}$ or $r=1$, then Theorem \ref{ratio} follows from Lemmas \ref{firstdn}--\ref{thirddn} just as in the argument following (\ref{finalcn}) for the case where $\GC = C_n$.
Hence it remains to handle the case where $\LC' = D_{n-r} \times A_{r-1}$ with $2\le r \le n$, $r \ne n-1$. We deal with this case in the next two lemmas.

\begin{lem}\label{fourthdn} Suppose $\LC' = D_{n-r} \times A_{r-1}$ with $r\ge 3$. Then $\a(\LC) \le \frac{1}{2}\left(1+\frac{\dim \LC}{\dim \GC}\right)$.
\end{lem}

\pf Let $u = u_1u_2$ be a unipotent element of $\LC$, where $u_1 \in D_{n-r}$, $u_2 \in A_{r-1}< D_r$.
We will show that
\begin{equation}\label{ineqdn}
\frac{\dim u_2^{A_{r-1}}}{\dim u_2^{D_r}} \le \frac{1}{2}\left(1+\frac{\dim \LC}{\dim \GC}\right).
\end{equation}
Given this, the lemma follows, since by Lemma \ref{thirddn} we have
\[
\frac{\dim u^\LC}{\dim u^\GC} \le \frac{\dim u_1^{D_{n-r}} + \dim u_2^{A_{r-1}}}{\dim u_1^\GC + \dim u_2^{D_r}},
\]
 and this is at most $\frac{1}{2}\left(1+\frac{\dim \LC}{\dim \GC}\right)$ by Lemma \ref{seconddn} and (\ref{ineqdn}).

It remains to establish (\ref{ineqdn}). Let $u_2 = \sum J_i^{a_i} \in A_{r-1} = SL_r$, so that $u_2$ has Jordan form
$\sum J_i^{2a_i}$ in $D_r$. Then
\[
\begin{array}{l}
\dim \CB_{GL_r}(u_2) = \sum ia_i^2+2\sum_{i<j}ia_ia_j,\\
\dim \CB_{D_r}(u_2) = 2\sum ia_i^2+4\sum_{i<j}ia_ia_j - \sum_{i\;odd}a_i.
\end{array}
\]
Write $s_2 :=  \sum_{i\;odd}a_i$. Then
$$\frac{\dim u_2^{A_{r-1}}}{\dim u_2^{D_r}}  = \frac{1}{2}(1+\frac{r-s_2}{\dim u_2^{D_r}}),$$
so to prove (\ref{ineqdn}) it suffices to show
\begin{equation}\label{ineqdn1}
\frac{r-s_2}{\dim u_2^{D_r}} \le \frac{\dim \LC}{\dim \GC} = \frac{2(n-r)^2-(n-r)+r^2}{2n^2-n}.
\end{equation}
It is straightforward to see that the right hand side of (\ref{ineqdn1}) is at least $\frac{1}{3}$, so (\ref{ineqdn1}) holds if
$\dim u_2^{D_r} \ge 3(r-s_2)$. The minimum value of $\dim u_2^{D_r}$ occurs when $u_2 = J_2+J_1^{r-2} \in SL_r$, in which case $\dim u_2^{D_r} = 4r-6$. This shows that (\ref{ineqdn1}) holds when $r\ge 6$.

It remains to establish (\ref{ineqdn1}) for $r = 3,4,5$. For $r=5$, the possibilities for $u_2 \in SL_r$ are as follows:
\[
\begin{array}{|c|c|c|c|c|c|c|}
\hline
\\[-13pt]
u_2 \in SL_5 & J_2+J_1^3 & J_2^2+J_1 & J_3+J_1^2 & J_3+J_2 & J_4+J_1 & J_5 \\
\hline
\\[-12pt]
\dim u_2^{D_5} & 14 & 20 & 26 & 28 & 32 & 36 \\
\hline
s_2 & 3 & 1 & 3 & 1 & 1 & 1 \\
\hline
\end{array}
\]
For all these possibilities (\ref{ineqdn1}) holds. The arguments for $r=3,4$ are similar. \hal

\begin{lem}\label{fifthdn} Suppose $\LC' = D_{n-2} \times A_{1}$. Then $\a(\LC) \le \frac{1}{2}\left(1+\frac{\dim \LC}{\dim \GC}\right)$.
\end{lem}

\pf Let $u = u_1u_2 \in \LC$ with $u_1 = \sum J_i^{n_i} \in D_{n-2}$ (so $\sum in_i = 2n-4$) and $u_2 \in A_1$. If $u_2=1$ then the conclusion follows from Lemma \ref{seconddn}, so assume $u_2 \ne 1$. Then $\dim u_2^{A_1} = 2$ and
$u = J_2^{n_2+2}+\sum_{i\ne 2}J_i^{n_i} \in \GC = D_n$. By \cite[3.1]{LS},
\[
\begin{array}{l}
\dim \CB_\LC(u) = \frac{1}{2}\sum in_i^2 + \sum_{i<j}in_in_j -\frac{1}{2}\sum_{i\;odd}n_i+2, \\
\dim \CB_\GC(u) = \frac{1}{2}\sum in_i^2 + \sum_{i<j}in_in_j -\frac{1}{2}\sum_{i\;odd}n_i +2n_1 + 4\sum_{i\ge 2}n_i+4.
\end{array}
\]
Define $s := \dim [V,u] = 2n-\sum n_i-2$. Then
\[
\dim u^\GC - \dim u^\LC = 4s+2n_1-8.
\]
Also \cite[3.4]{lish99} gives $\dim u^\GC \le \frac{1}{2}s(4n-s+1)$. Hence we see that the desired inequality
$ \frac{\dim u^\LC}{\dim u^\GC} \le \frac{1}{2}\left(1+\frac{\dim \LC}{\dim \GC}\right)$ is equivalent to the following
\begin{equation}\label{lastdn}
 \frac{1}{2}s(4n-s+1) \le \frac{(2n^2-n)(4s+2n_1-8)}{4n-7}.
\end{equation}
Now $2n-4 = \sum in_i \ge 2\sum_{i\ge 2}n_i$, and hence
\[
s = 2n-2 - n_1 -\sum_{i\ge 2}n_i \ge 2n-2-n_1 -(n-2) = n-n_1.
\]
It follows that $4s+2n_1 \ge 2n+2s$, and hence (\ref{lastdn}) holds provided
$$\frac{1}{2}s(4n-s+1) \le \frac{(2n^2-n)(2n+2s-8)}{4n-7},$$
which is true for all $s$ when $n\ge 5$.
Finally, when $n=4$ the conclusion of the lemma is easily checked directly. This completes the proof. \hal

\vspace{4mm}
This completes the proof of Theorem \ref{ratio}.

\subsection{Exceptional groups: Proof of Theorem \ref{alphaexcep}}

Let $\GC$ be a simple algebraic group of exceptional type in good characteristic. In principle one can check Theorem \ref{alphaexcep} by going through all possible Levi subgroups $\LC$ of $\GC$, in each case listing all the unipotent class representatives $u$  in $\LC$ and using Theorem 3.1 and Tables 22.1.1--5 of \cite{LS} to write down the dimensions of $u^\LC$ and $u^\GC$. In fact, this is precisely what we do for the Levi subgroups listed in Table \ref{extab}, and for the remaining ones (labelled ``rest" in Table \ref{extab}) we need a short argument.

We will give the proof of Theorem \ref{alphaexcep} just for $\GC = E_7$ and leave the other entirely similar cases to the reader. First suppose that the Levi subgroup $\LC$ is one of those listed for $\GC=E_7$ in Table \ref{extab}. In each case we adopt the above procedure of listing unipotent representatives $u$  in $\LC$ and calculating $\dim u^\LC$ and $\dim u^\GC$. We illustrate below with the case $\LC' = D_6$, listing in the first row the Jordan form of $u$ on the 12-dimensional module for $\LC'$ and in the second row the class of $u$ in $\GC$ as in \cite[Table 22.1.2]{LS}:

{\small
\[
\begin{array}{|c|c|c|c|c|c|c|c|}
\hline
\\[-11pt]
u \hbox{ in }D_6 & (2^2,1^8) & (3,1^9) & (2^4,1^4) & (2^6) & (2^6)' & (3,2^2,1^5) & (3,2^4,1)\\
u \hbox{ in }E_7 & A_1 & A_1^2 & A_1^2 & (A_1^3)^{(1)} & (A_1^3)^{(2)} & (A_1^3)^{(2)} & A_1^4  \\
\dim u^\LC & 18 & 20 & 28 & 30 & 30 & 32 & 36  \\
\dim u^\GC & 34 & 52 & 52 & 54 & 64 & 64 & 70  \\
\hline
\end{array}
\]
\[
\begin{array}{|c|c|c|c|c|c|c|c|}
\hline
\\[-11pt]
 (3,2,1^6) & (3^2,2^2,1^2) & (3^3,1^3) & (3^4) & (4^2,1^4) &  (5,1^7) & (4^2,2^2) & (4^2,2^2)' \\
  A_2 & A_2A_1 & A_2A_1^2 & A_2^2 & A_3 & A_3 & (A_3A_1)^{(1)} & (A_3A_1)^{(2)} \\
  34 & 40 & 42 & 44 & 44 & 36 & 46 & 46 \\
  66 & 76 & 82 & 84 & 84 & 84 & 86 & 92 \\
\hline
\end{array}
\]
\[
\begin{array}{|c|c|c|c|c|c|c|c|}
\hline
\\[-11pt]
(5,2^2,1^3) & (4^2,3,1) & (5,3^2,1) & (5^2,1^2) & (6^2) & (6^2)' & (7,1^5) & (5,3,1^4) \\
 (A_3A_1)^{(1)} & A_3A_1^2 & A_3A_2 & A_4 & (A_5)^{(1)} & (A_5)^{(2)} & D_4 & D_4(a_1) \\
44 & 48 & 50 & 52 & 54 & 54 & 48 & 46 \\
86 & 94 & 98 & 100 & 102 & 108 & 96 & 94 \\
\hline
\end{array}
\]
\[
\begin{array}{|c|c|c|c|c|c|c|}
\hline
\\[-11pt]
(7,2^2,1) & (5,3,2^2) & (9,1^3) & (7,3,1^2) & (11,1) & (9,3,) & (7,5) \\
D_4A_1 & D_4(a_1)A_1 & D_5 & D_5(a_1) & D_6 & D_6(a_1) & D_6(a_2) \\
52 & 48 & 56 & 54 & 60 & 58 & 56 \\
102 & 96 & 112 & 106 & 118 & 114 &  110 \\
\hline
\end{array}
\]
}

To compute the information in the tables, we list the possible Jordan forms for unipotent elements $u$ of $D_6$, and in each case find a Levi subgroup of $D_6$ in which $u$ is contained as a regular element; this Levi subgroup then gives the label of $u$ as an element of $E_7$ in Table 22.1.2 of \cite{LS}. For cases where all the Jordan blocks have even size -- namely the Jordan forms $(2^6)$, $(4^2,2^2)$ and $(6^2)$ -- there are two $D_6$-classes (see \cite[3.11]{LS}), and the corresponding $E_7$-classes can be worked out by computing the dimension of $\CB_{L(E_7)}(u)$ using the restriction $L(E_7)\downarrow D_6$ (see \cite[11.8]{LS}), where $L(E_7)$ denotes the Lie algebraa of $\GC = E_7$.

Inspecting the tables above, we see that the maximum value of $\frac{\dim u^\LC}{\dim u^\GC}$ is equal to $\frac{30}{54}$, and is attained when $u$ has Jordan form $(2^6)$ in $D_6$  and is in the class $(A_1^3)^{(1)}$ of $E_7$. Hence for $\LC'=D_6$ we have $\a(\LC) = \frac{5}{9}$, as in Table \ref{extab} of Theorem \ref{alphaexcep}.

Now suppose $\LC$ is not one of the Levi subgroups listed for $\GC=E_7$ in Table \ref{extab} (i.e. $\LC$ does not have a factor $E_6$, $D_r$ or $A_r$ ($r\ge 3$). Then $\dim \LC \le \dim A_2A_2A_1T_2 = 21$.

 Let $u$ be a nontrivial unipotent element in $\LC$, and assume for a contradiction that
\[
\frac{\dim u^\LC}{\dim u^\GC} > \frac{1}{6}.
\]
Now $\dim u^\LC \le \dim \LC-7\le 14$, and hence $\dim u^\GC < 84$. It then follows from Table 22.1.2 of \cite{LS} that $u$ is in one of the following $E_7$-classes :
\[
A_1, A_1^2, (A_1^3)^{(1)}, (A_1^3)^{(2)}, A_1^4, A_2, A_2A_1, A_2A_1^2.
\]
For these classes the maximum possible value of $\dim u^\LC$ occurs for $\LC' = A_2A_2A_1$ or $A_1^4$, and is as follows:\\
{\small
\[
\begin{array}{|c|c|c|c|c|c|c|c|c|}
\hline
\\[-11pt]
u & A_1 &  A_1^2 & (A_1^3)^{(1)} & (A_1^3)^{(2)} & A_1^4 & A_2 & A_2A_1 & A_2A_1^2 \\
\hline
\\[-11pt]
\hbox{ max dim }u^\LC, & 4 & 8 & - & 10 & - & 6 & 10 & 12 \\
\LC' = A_2A_2A_1 &&&&&&&& \\
\hline
\\[-11pt]
\hbox{ max dim }u^\LC,& 2 & 4 & 6 & 6 & 8 & - & -& - \\
\LC' = A_1^4 &&&&&&&& \\
\hline
\\[-11pt]
\dim u^\GC & 34 & 52 & 54 & 64 & 70 & 66 & 76 & 82 \\
\hline
\end{array}
\]
}
In all cases we see that $\frac{\dim u^\LC}{\dim u^\GC} < \frac{1}{6}$, which is a contradiction. This completes the proof of Theorem \ref{alphaexcep} for $\GC = E_7$. \hal

\subsection{Proof of Corollary \ref{keycor} and Theorem \ref{supp}}
The proof of Corollary \ref{keycor} is immediate, since
\[
1-\frac{1}{2}\frac{\dim y^{\GC}}{\dim \GC} = \frac{1}{2}\left(1+ \frac{\dim \CB_\GC(y)}{\dim \GC}\right) =
 \frac{1}{2}\left(1+ \frac{\dim \LC}{\dim \GC}\right),
\]
and the right hand side above is at least $\a(\LC)$ by Theorems \ref{ratio} and \ref{alphaexcep}.

\hal
\medskip

To prove Theorem \ref{supp}, note that $\c(1) \geq q^r/3$ by \cite{LSe}.
Also, $\CB_G(g) \leq \CB_G(y) = L$, so by Theorem \ref{main1} and the inequality
$\a(L) \le 1-\frac{1}{2}\frac{\dim y^{\GC}}{\dim \GC}$ obtained above, we have
$$|\c(g)| \le f(r) \, \c(1)^{1-\frac{1}{2}\frac{\dim y^{\GC}}{\dim \GC}}.$$
Hence it suffices to prove that $\g r \geq cs$, where $\g := (\dim y^\GC)/(2\dim \GC)$ and $s: = \supp(y)$.
Define $a:=1$ if $\GC = SL_n$ and $a:=1/2$ otherwise.

Lemma 3.4 of \cite{lish99} relates the support of elements of prime order in $G$ with the
size of their conjugacy class. The proof of this lemma only uses the fact that these elements
are semisimple or unipotent. Since $y \in \ZB(L)$ is semisimple, the lemma applies and shows
in particular that $|y^G| \ge c' q^{ans}$, where $c' >0$ is an absolute constant.
This implies that $\dim y^\GC \ge ans$, and so
\[
\g r \geq \frac{ans}{2\dim \GC}r = \frac{anr}{2\dim \GC}s = c s,
\]
as needed.
\hal
\medskip

\subsection{Bounds for $GL_n$: proof of Theorem \ref{GL} and Corollary \ref{FL}}

Let $\KK$ an algebraically closed field of characteristic $p$,
and let $\LC = GL_{n_1}(\KK) \times \cdots \times GL_{n_m}(\KK)$, so that the Levi subgroup
$L$ in Theorem \ref{GL} can be viewed as $\LC^F$ for a suitable Frobenius endomorphism $F$.
Fix $n$ pairwise distinct elements $\la_1, \ldots,\la_n \in \KK^\times$. The statements follow from
Theorem \ref{main1} if $n_{i_0} = 1$, so we will assume that $n_{i_0} \geq 2$.

Any unipotent element $u \in L$ can be written as $\diag(u_1, \ldots, u_m)$, where $u_i \in \LC_i := GL_{n_i}(\KK)$ is
unipotent. Write $u_i = J_{b_{i1}} \oplus \ldots \oplus J_{b_{ir_i}}$ for a partition
$\nu_i := (b_{i1} \geq b_{i2} \geq \ldots \geq b_{ir_i} \geq 1)$ of $n_i$, and define
$$s_i := \diag(\la_1, \la_2,\ldots,\la_{b_{i1}},\la_1, \la_2,\ldots,\la_{b_{i2}}, \ldots, \la_1, \la_2,\ldots,\la_{b_{ir_i}}) \in \LC_i$$
Note that if $a_{ij}$ is the multiplicity of $\la_j$ as an eigenvalue of $s_i$, then $(a_{i1} \geq a_{i2} \geq \ldots \geq a_{in} \geq 0)$ is the partition of
$n_i$ conjugate to $\nu_i$. Now Lemma \ref{uss} shows that
$$\dim u_i^{\LC_i} = \dim s_i^{\LC_i} = n_i^2 - \sum^n_{j=1}a_{ij}^2.$$
Similarly, setting $s := \diag(s_1, \ldots, s_m) \in \LC$, we then get
$$\dim u^\GC = \dim s^\GC = n^2 - \sum^n_{j=1}(\sum^m_{i=1}a_{ij})^2.$$
Note that $u \neq 1$ precisely when $\max_{1 \leq i \leq m}a_{i2} > 0$. Thus
$\a(L) = \beta(n_1, \ldots, n_m)$. Now
Theorem \ref{GL} follows immediately from Theorem \ref{main1} and Theorem \ref{GL2}(i) below,

\smallskip
Note by the Cauchy-Schwartz inequality that
$(\sum_{i = 1}^m a_{ij})^2 \le m \sum_{i=1}^m a_{ij}^2$ for each $j$, with equality attained exactly when
$a_{1j} = a_{2j} = \ldots = a_{mj}$.
Setting $\Delta := \sum_{i, j} a_{ij}^2$, we have
\begin{equation}\label{4.3}
\beta(n_1, \ldots , n_m) \le \max_{\Delta} \frac{ (\sum_{i=1}^m n_i^2) - \Delta }{ n^2 - m\Delta}.
\end{equation}
Now suppose $n_i = n/m$ for $i = 1, \ldots , m$. Then $\sum_{i=1}^m n_i^2/n^2 = 1/m$, and so
\eqref{4.3} implies that $\beta(n_1, \ldots,n_m) \leq 1/m$. In fact equality holds if we choose $a_{1j} = a_{2j} = \ldots = a_{mj}$ for all $j$.
Thus Corollary \ref{FL} follows.
\hal

In what follows, for any partitions $\a = (a_1 \geq a_2 \geq \ldots \geq a_n \geq 0) \vdash A$ and $\b = (b_1 \geq b_2 \geq \ldots \geq b_n \geq 0) \vdash B$
of $A, B \geq 1$ we set
\begin{equation}\label{gl21-1}
  g(\a) := A^2-\sum^n_{i=1}a_i^2,~~h(\a) := \frac{g(\a)}{A},~~\a+\b := (a_1+b_1,a_2+b_2, \ldots,a_n+b_n) \vdash (A+B).
\end{equation}

\begin{lem}\label{GL21}
Let $\a = (a_1 \geq a_2 \geq \ldots \geq a_n \geq 0) \vdash A$ and $\b = (b_1 \geq b_2 \geq \ldots \geq b_n \geq 0) \vdash B$ be two
partitions of $A, B \geq 1$. Then $h(\a) + h(\b) \leq h(\a+\b)$.
\end{lem}

\pf
We need to show that
$$\frac{A^2-\sum^n_{i=1}a_i^2}{A}+\frac{B^2-\sum^n_{i=1}b_i^2}{B} \leq \frac{(A+B)^2-\sum^n_{i=1}(a_i+b_i)^2}{A+B},$$
equivalently, $\G \geq 0$, where
$$\begin{aligned}\G & := AB\sum_{i \neq j}(a_i+b_i)(a_j+b_j)-(A+B)(B\sum_{i \neq j}a_ia_j+A\sum_{i \neq j}b_ib_j)\\
    & = AB\sum_{i \neq j}(a_ib_j+a_jb_i)-(B^2\sum_{i \neq j}a_ia_j+A^2\sum_{i \neq j}b_ib_j)\\
    & = AB(2AB-2\sum_ia_ib_i)-B^2(A^2-\sum_ia_i^2)-A^2(B^2-\sum_ib_i^2)\\
    & = B^2\sum_ia_i^2+A^2\sum_ib_i^2-2AB\sum_ia_ib_i.
    \end{aligned}$$
By the Cauchy-Schwarz inequality,
$$2AB\sum_ia_ib_i \leq 2\cdot B(\sum_ia_i^2)^{1/2} \cdot A (\sum_ib_i^2)^{1/2} \leq  B^2\sum_ia_i^2+A^2\sum_ib_i^2,$$
and the claim follows.
\hal

\begin{thm}\label{GL2}
In the notation of Theorem \ref{GL}, assume that
$$n_1 = n_2 = \ldots = n_t > n_{t+1} \geq \ldots \geq n_m \geq 1.$$
Then the following statements hold.
\begin{enumerate}[\rm(i)]
\item $n_1/n \geq \beta(n_1, \ldots,n_m) \geq (n_1-1)/(n-t).$
\item If $m=2$, then $\beta(n_1,n_2) = (n_1-1)/(n-t)$. Moreover, if
$$1 \neq u = \diag(u_1,u_2) \in GL_{n_1}(q) \times GL_{n_2}(q) = L$$
is a unipotent element, then
$(\dim u^\LC)/(\dim u^\GC) = \a(L)$ precisely when one of the following conditions holds.
\begin{enumerate}[\rm(a)]
\item $n_1 = n_2$, and $u_1$ and $u_2$ have the same Jordan canonical form.
\item $n_1 > n_2$, $u_1$ is a transvection and $u_2 = 1$.
\item $n_1 = n_2+1$, and the sizes of Jordan blocks for $u_1$ and $u_2$ are
$$c_1 \geq \ldots \geq c_{j-1} \geq c_j \geq c_{j+1} \geq \ldots \geq c_s,~~
    c_1 \geq \ldots \geq c_{j-1} \geq c_j-1 \geq c_{j+1} \geq \ldots \geq c_s,$$
respectively.
\end{enumerate}
\end{enumerate}
\end{thm}

\pf
(i) To prove the lower bound for $\b(n_1,\ldots,n_m)$, we choose $(a_{i1},\ldots,a_{in})$ to be $(n_1-1,1,0,\ldots,0)$ if $1\leq i \leq t$ and $(n_i,0,\ldots,0)$ otherwise. To prove the upper bound, for $1 \leq i \leq m$ consider the partition $\a_i:=(a_{i1},a_{i2}, \ldots,a_{in}) \vdash n_i$. By Lemma
\ref{GL21} we have
\begin{equation}\label{gl21-2}
  \sum^m_{i=1}\frac{g(\a_i)}{n_i} = h(\a_1)+h(\a_2) + \ldots + h(\a_m) \leq h(\sum^m_{i=1}\a_i) = \frac{g(\sum^m_{i=1}\a_i)}{n}
\end{equation}
in the notation of \eqref{gl21-1}. Note that the condition $\max_{1 \leq i \leq m} a_{i2} > 0$ ensures that $g(\sum^m_{i=1}\a_i) > 0$.
Since $n_1 = \max_{1 \leq i \leq m} n_i$,
\eqref{gl21-2} now implies that
$$\beta(n_1,n_2, \ldots,n_m) = \frac{\sum^m_{i=1}g(\a_i)}{g(\sum^m_{i=1}\a_i)} \leq
   \frac{n_1(\sum^m_{i=1}g(\a_i)/n_i)}{g(\sum^m_{i=1}\a_i)} \leq \frac{n_1}{n}.$$

\smallskip
(ii) We may assume $A := n_1 > B:= n_2$ by  Corollary \ref{FL} and its proof. To ease the notation, also write
$$(a_{11},a_{12}, \ldots,a_{in}) = (a_1, a_2, \ldots,a_n),~~(a_{21},a_{22}, \ldots,a_{2n}) = (b_1, b_2, \ldots,b_n).$$
Then we need to show that
$$\frac{\sum_{1 \leq i \neq j \leq n}(a_ia_j+b_ib_j)}{\sum_{1 \leq i \neq j \leq n}(a_i+b_i)(a_j+b_j)} \leq \frac{A-1}{A+B-1},$$
equivalently, $\Sigma \geq 0$, where
$$\begin{aligned}\Sigma & := (A-1)\sum_{i \neq j}(a_i+b_i)(a_j+b_j)-(A+B-1)\sum_{i \neq j}(a_ia_j+b_ib_j)\\
    & = (A-1)\sum_{i \neq j}(a_ib_j+a_jb_i)-B\sum_{i \neq j}(a_ia_j+b_ib_j)\\
    & = (A-1)(2AB-2\sum_ia_ib_i)-B(A^2+B^2-\sum_ia_i^2-\sum_ib_i^2)\\
    & = B(A^2-2A-B^2)+B(\sum_ia_i^2+\sum_ib_i^2)-2(A-1)\sum_ia_ib_i\\
    & = B(A^2-2A-B^2)+B\sum_i(a_i-b_i)^2-2(A-1-B)\sum_ia_ib_i\\
    & = (A-1-B)((A-1+B)B-2\sum_ia_ib_i)+B(\sum_i(a_i-b_i)^2-1)\\
    & = (A-1-B)(\sum_i(A-1+B-2a_i)b_i)+B(\sum_i(a_i-b_i)^2-1).
    \end{aligned}$$
Note that the condition $A \geq B+1$ implies that $\sum_i(a_i-b_i)^2 \geq 1$, with equality attained exactly when
\begin{equation}\label{GL2-1}
  A=B+1,~~(a_1,a_2, \ldots,a_n) = (b_1, \ldots ,b_{i-1},b_i+1,b_{i+1}, \ldots,b_n).
\end{equation}
First we consider the case when $A-1+B \geq 2a_i$ for all $i$. As $B \geq 1$, we see that $\Sigma \geq 0$, with
equality attained exactly when \eqref{GL2-1} holds, which means the corresponding unipotent element satisfies (c).

Suppose now that $A-1+B \leq 2a_i-1$ for some $i$. Then $a_i \geq \sum_{j \neq i}a_j+B$. As $B \geq 1$ and
$a_1 \geq a_2 \geq \ldots \geq a_n$, this can happen only for one index $i$, and this index $i$ is $1$, and so
\begin{equation}\label{GL2-2}
  a_1 \geq A'+B,
\end{equation}
where $A' := \sum_{j \geq 2}a_j$. In particular, $A-1+B \geq 2a_j$ for all $j \geq 2$. Now by \eqref{GL2-2} we have
$$\begin{aligned} \Sigma' := & (A-1-B)((A-1+b_1)b_1-2a_1b_1)+b_1(a_1-b_1)^2\\
    = &~b_1((A-1-B)(b_1-1+A'-a_1)+(a_1-b_1)^2)\\
    = &~b_1((a_1-b_1)(A'-1)+(A'-1-B')(b_1-1+A'-a_1))\\
    = &~b_1((a_1-b_1)B'+(A'-1)^2-B'(A'-1))\\
    = &~b_1((a_1-A'-b_1+1)B'+(A'-1)^2) \geq 0,
    \end{aligned}$$
with equality exactly when
\begin{equation}\label{GL2-3}
  B' := \sum_{j \geq 2}b_j = 0,~~A' =1.
\end{equation}
It follows that, if $\sum_{j \geq 2}(a_j-b_j)^2 \geq 1$, then $\Sigma \geq 0$, with equality exactly when \eqref{GL2-3} holds,
which means the corresponding unipotent element satisfies (b).

Assume finally that $\sum_{j \geq 2}(a_j-b_j)^2 \leq 0$. Then $a_j=b_j$ for all $j \geq 2$ and $A'=B'$. As $\max(a_2,b_2) > 0$,
we must have $A'=B' \geq 1$, and so by \eqref{GL2-2}
$$\Sigma' = b_1((a_1-b_1-1)A'+1) \geq 2b_1,$$
yielding $\Sigma \geq b_1 \geq 1$.
\hal

We note that Theorem \ref{GL2}(i) was inspired by some correspondence with M. Fraczyk who is studying the situation in
Theorem \ref{GL} using different methods.

\section{Random walks}

In this section we prove Theorems \ref{mix}--\ref{mix-gluni} concerning random walks and covering numbers.

\vspace{2mm}
\noindent {\bf Proof of Theorem \ref{mix}}

Suppose $\GC$ is a  simple algebraic group of rank $r$ in good characteristic,
and $G = G(q) = \GC^F$ is a finite quasisimple group over $\F_q$.
Let $y \in G$ be such that $\CB_G(y) \le  L = \LC^F$ for a split
Levi subgroup $\LC$ of $\GC$. Write $C = y^G$, and let $h$ be the Coxeter number of $\GC$.

For a real number $s$, define
\begin{equation}\label{zetadef}
\zeta^G(s) = \sum_{\c \in \Irr(G)} \c(1)^{-s}.
\end{equation}
We will need the following result, which is \cite[Theorem 1.1]{chardeg}.

\begin{lem}\label{zetabd}
If $s>\frac{2}{h}$, then $\zeta^G(s) \rightarrow 1$ as $q\rightarrow \infty$.
\end{lem}

We first prove part I(a) of Theorem \ref{mix} together with the first statement of part (II) (the $C^6=G$ statement). We will prove the mixing time assertions later.

Let $t$ be a positive integer. By a well-known result (see \cite[Chapter 1, 10.1]{AH}),
for $g \in G$  the number of ways of writing $g$ as a product of $t$ conjugates of $y$ is
\[
N(g) = \frac{|C|^t}{|G|} \sum_{\c \in \Irr(G)} \frac{\c(y)^t \c(g^{-1})} {\c(1)^{t-1}}.
\]
Define $P^t(g) = \frac{N(g)}{|C|^t}$, the probability that a random product of $t$ conjugates of $y$ is equal to $g$,
and let $U(g) = \frac{1}{|G|}$, the uniform probability distribution on $G$.
Then
\begin{equation}\label{one}
|P^t(g)- U(g)| \le  \frac{1}{|G|}\sum_{\c(1)> 1} \left(\frac{|\c(y)|}{\c(1)}\right)^t\c(1)^2.
\end{equation}
Define
\[
||P^t-U||_\infty = |G|\,\hbox{max}_{g\in G}|P^t(g)-U(g)|.
\]
Write $\a = \a(\LC)$. Then Theorem \ref{main1} gives $\frac{|\c(y)|}{\c(1)} \le f(r)\c(1)^{\a-1}$, and so (\ref{one}) implies
\[
\begin{array}{ll}
||P^t-U||_\infty  & \le f(r)^t \sum_{\c(1)> 1} \c(1)^{t(\a-1)+2} \\
                    &  = f(r)^t\left(\zeta^G\left(t(1-\a)-2\right)-1\right).
\end{array}
\]
By Lemma \ref{zetabd}, $\zeta^G(t(1-\a)-2)-1 \rightarrow 0$ as $q\rightarrow \infty$ provided
\begin{equation}\label{ineqal}
t(1-\a)-2 > \frac{2}{h}.
\end{equation}

If $\GC$ is of exceptional type $G_2,F_4,E_6,E_7$ or $E_8$, then $\frac{2}{h}$ is $\frac{1}{3}$, $\frac{1}{6}$,
$\frac{1}{6}$, $\frac{1}{9}$ or $\frac{1}{15}$ respectively, and Theorem \ref{alphaexcep} shows that (\ref{ineqal}) holds in all cases, provided $t \ge 6$. This proves the first statement of Theorem \ref{mix}(II).

Now suppose $\GC$ is of classical type. Then $\a \le \frac{1}{2}(1+\frac{\dim \LC}{\dim \GC})$ by Theorem \ref{ratio}. This implies that (\ref{ineqal}) holds provided $t > (4+\frac{4}{h})\frac{\dim \GC}{\dim \GC - \dim \LC}$, proving Theorem \ref{mix}(I)(a).

\vspace{2mm}
We now prove the assertions on mixing times in Theorem \ref{mix}. For these  we use the Diaconis--Shashahani bound
\cite{DS}:
\begin{equation}\label{upperbd}
(||P^t-U||_1)^2 \le \sum_{\c \in \Irr(G), \c\ne 1} \left(\frac{|\c(y)|}{\c(1)}\right)^{2t}  \c(1)^2.
\end{equation}
As above, Theorem \ref{main1} shows that the right hand side of (\ref{upperbd}) is less than \linebreak     $f(r)^{2t}\left(\zeta^G(2t(1-\a)-2)-1\right)$, and hence tends to 0 as $q\rightarrow \infty$ provided $2t(1-\a)-2 > \frac{2}{h}$. Using Theorems \ref{alphaexcep} and \ref{ratio}, we now see as before that this inequality holds provided $t \ge 3$ when $\GC$ is of exceptional type, and provided $t > (2+\frac{2}{h})\frac{\dim \GC}{\dim \GC - \dim \LC}$ when $\GC$ is classical. This proves the mixing time assertions, completing the proof of Theorem \ref{mix}. \hal

\vspace{4mm}
\noindent {\bf Proof of Theorem \ref{mix-gluni}}

This is very similar to the previous proof, using Corollary \ref{sl-uni} instead of Theorem \ref{main1}.
Let $G = SL_n(q)$ and let $u \in G$ be a non-identity unipotent element.
Let $t \in \N$, and for $g\in G$ let $P^t(g)$ be the probability that a random product of $t$ conjugates of $u$ is equal to $g$, and $U(g) = \frac{1}{|G|}$. As in (\ref{one}),
\[
|P^t(g)- U(g)| \le  \frac{1}{|G|}\sum_{\c \in \Irr(G),\c(1)> 1} \left(\frac{|\c(u)|}{\c(1)}\right)^t\c(1)^2.
\]
By Corollary \ref{sl-uni}, $\frac{|\c(u)|}{\c(1)} \le g(n)\,\c(1)^{-\frac{1}{n-1}}$ for $\c \in \Irr(G)$, and hence
\[
\begin{array}{ll}
||P^t-U||_\infty  & \le g(n)^t \sum_{\c(1)> 1} \c(1)^{-\frac{t}{n-1}+2} \\
                    &  = g(n)^t\left(\zeta^G(\frac{t}{n-1}-2)-1\right).
\end{array}
\]
By Lemma \ref{zetabd}, $\zeta^G(\frac{t}{n-1}-2)-1 \rightarrow 0$ as $q\rightarrow \infty$ provided
$\frac{t}{n-1}-2 > \frac{2}{n}$, which holds provided $t > 2n$. This proves part (i) of Theorem \ref{mix-gluni}. Part (ii) is proved in the same way, using the bound (\ref{upperbd}). \hal

\vspace{4mm}
\noindent {\bf Corollaries \ref{linearbd} and \ref{linearbdsl}}

Corollary \ref{linearbd} follows immediately from Theorem \ref{mix}(I)(b). Corollary \ref{linearbdsl} is proved exactly as above, using Theorem \ref{main1c}.

\vspace{4mm}
Next, we use some well-known observations to justify the remarks made after the statement of Theorem \ref{mix}.

\begin{lem}\label{mix-subset}
\begin{enumerate}[\rm(i)]
\item Let $G$ be a finite group, and let $S$ be a generating subset of
$G$ that satisfies $|S^N| < |G|(1-1/e)$
for some integer $N \geq 1$. Then the mixing time $T(G,S)$ of the random walk
on the Cayley graph corresponding to $S$ is at least $N+1$.
\item Let $G = SL_n(q)$ with $n \geq 2$ and $S = y^G$ with
$y = \diag(\mu I_{n-1},\la)$, where $\mu, \la \in \F_q^\times$ and $\mu \neq \la$.
Then $T(G,y) \geq n$.
\end{enumerate}
\end{lem}

\pf
(i) Define $P(g)$ to be $1/|S|$ if $g \in S$ and $0$ otherwise, and
let $U(g) = 1/|G|$ for all $g \in G$. Consider any $1 \leq k \leq N$. Note that
$|S^k| \leq |S^{k+1}|$ and so $|S^k| \leq |S^N|$, whence
$$||P^k-U||_1 \geq \sum_{g \in G \smallsetminus S^k}|P^k(g)-U(g)|
  = \sum_{g \in G \smallsetminus S^k}|U(g)|
  \geq \frac{|G \smallsetminus S^k|}{|G|} > 1/e.$$
It follows that $T(G,S) \geq N+1$.

(ii) Note that $(y^G)^{n-1}$ is contained in $X$, the set of
elements $x \in G$ that have eigenvalue $\mu$ on $V = \F_q^n$. Now if we fix
$0 \neq v \in V$ and let $Y := \{ x \in G \mid x(v) = \mu v\}$,
then it is easy to see that $|\NB_G(Y)|/|Y| \geq q-1 \geq 2$. Hence,
$$|X| = |\cup_{g \in G}gYg^{-1}| \leq |Y| \cdot [G:\NB_G(Y)] \leq |G|/2.$$
Now we can apply (i) to $S := y^G$.
\hal

We conclude with a proof of our last theorem, connecting the mixing times of random walks
on classical groups with the support of certain elements.

\vspace{2mm}
\noindent {\bf Proof of Theorem \ref{mix-supp}}

Set $s:= \supp(y)$. Then $\CB_G(g) \le \CB_G(y) = L$.
Theorem \ref{mix} I(b) gives
$$T(G,g) \le \lceil (2+\frac{2}{h})\frac{\dim \GC}{\dim \GC - \dim \LC}\rceil$$
for large $q$. Now, $\dim \GC - \dim \LC = \dim \GC - \dim \CB_{\GC}(y) = \dim y^{\GC} \ge ans$ as shown in the proof of
Theorem \ref{supp}. This yields
\[
T(G,g) \le \lceil (2+\frac{2}{h})\frac{\dim \GC}{ans}\rceil.
\]
Let $c = c(\GC)$ be as in Theorem \ref{supp}. Then we have $\frac{\dim \GC}{an} = \frac{r}{c} = r'$.
We obtain
\[
T(G,g) \le \lceil (2+\frac{2}{h})r' / s \rceil,
\]
proving the first assertion.

It remains to prove the lower bound on $T(G,y)$.
By (\ref{time}) we have
\[
T(G,y) \ge \frac{\log |G| + \log (1 - e^{-1})}{\log |y^{G}|} \gtrsim_{|G|} \frac{\dim \GC}{\dim y^{\GC}}.
\]
It follows from \cite[3.4]{lish99} and its proof that, for $y$ semisimple, we have $|y^G| \le 2ans$.
Hence $\dim y^{\GC} \le 2ans$, which, combined with the inequality above, implies
\[
T(G,y) \gtrsim_{|G|} \frac{\dim \GC}{2ans} \ge \frac{\dim \GC}{2an}/s = \frac{1}{2} r' / s,
\]
as required.
\hal

\end{document}